\documentclass[aos,MSNbibl,seceqn,nameyear,dvips]{arximspdf}

\doi{10.1214/12-AOS995} 
\volume{40}
\issue{2}
\pubyear{2012}
\firstpage{1132}
\lastpage{1170}

\makeatletter
\newcommand{\eqref}[1]{(\ref{#1})}
\newcommand{\E}{\mathbb{E}}
\newcommand{\IR}{\mathbb{R}}
\newcommand{\Var}{\operatorname{Var}}
\newcommand{\argmin}{\operatorname{argmin}\limits}
\newcommand{\ASF}{\operatorname{ASF}}
\newcommand{\diag}{\operatorname{diag}}

\newproclaim{remark}{Remark}
\newproclaim{assumption}{Assumption}
\newtheorem{theorem}{Theorem}
\newtheorem{proposition}{Proposition}
\newtheorem{corollary}{Corollary}
\newtheorem{lemma}{Lemma}

\makeatother

\begin{document}
\begin{frontmatter}

\title{Nonparametric regression with nonparametrically generated covariates}
\runtitle{Nonparametrically generated covariates}

\begin{aug}
\author[A]{\fnms{Enno} \snm{Mammen}\ead[label=e1]{emammen@rumms.uni-mannheim.de}},
\author[B]{\fnms{Christoph} \snm{Rothe}\corref{}\ead[label=e2]{rothe@cict.fr}\ead[label=u1,url]{http://www.christophrothe.net}}
\and
\author[C]{\fnms{Melanie} \snm{Schienle}\ead[label=e3]{melanie.schienle@wiwi.hu-berlin.de}}

\runauthor{E. Mammen, C. Rothe and M. Schienle}
\affiliation{University of Mannheim, Toulouse School of Economics and
Humboldt~University~Berlin}

\address[A]{E. Mammen\\
Department of Economics\\ University of Mannheim\\ D-68131 Mannheim\\
Germany\\
\printead{e1}} 
\address[B]{C. Rothe\\ Toulouse School of Economics\\ 21 Allee de
Brienne\\ F-31000 Toulouse\\ France\\
\printead{e2}\\
\printead{u1}}
\address[C]{M. Schienle\\ School of Business and Economics\\ Humboldt
University Berlin\\ Spandauer Str. 1\\ D-10178 Berlin\\
Germany\\ \printead{e3}}
\end{aug}

\received{\smonth{11} \syear{2011}}

%
\begin{abstract}
We analyze the statistical properties of nonparametric regression
estimators using covariates which are not directly observable,
but have be estimated from data in a preliminary step.
These so-called generated covariates appear in numerous applications,
including two-stage nonparametric regression, estimation of simultaneous
equation models or censored regression models. Yet so far there seems to
be no general theory for their impact on the final estimator's
statistical properties.
Our paper provides such results. We derive a stochastic expansion that
characterizes the influence of the generation step on the final estimator,
and use it to derive rates of consistency and asymptotic distributions
accounting for the presence of generated covariates.
\end{abstract}

%


\end{frontmatter}

\section{Introduction}
A wide range of statistical applications requires nonparametric
estimation of a regression function when some of the covariates
are not directly observed, but have themselves only been
estimated in a~(possibly nonparametric) preliminary step.
Examples include triangular simultaneous equation models
[e.g., \citet{newey1999nonparametric}, \citet{blundell2004endogeneity}, \citet{imbens2009identification}], sample selection models
[\citet{das2003nonparametric}], treatment effect models
[\citet{heckman1998matching}, \citet{heckman2005structural}], censored
regression models [\citet{linton2002nonparametric}], generalized
Roy models [\citet{dHautfoeuille2009Roy}], stochastic volatility
models [\citet{kristensen2009}] and GARCH-in-Mean models
[\citet{conrad2009garch}], amongst many others. In
contrast to fully parametric settings [\citet{pagan1984econometric}],
there seems to be no general theoretical results on how to derive
the statistical properties of such
nonparametric two-step estimators. Instead, most available
results in the literature typically exploit peculiarities of a
specific model, and can thus not easily be transferred to other
applications.

In this paper, we study the statistical properties of a
nonparametric estimator $\hat m_{LL}$ of a conditional mean
function $m_0(x) = \E(Y|r_0(S) = x)$ when the function $r_0$ is
unknown, but can be estimated from data. While we are specific
about estimating $m_0$ by local linear regression
[\citet{fan1996local}] to simplify technical arguments, we
neither require the generated regressors $\hat{R} = \hat{r}(S)$
to emerge from a specific type of model, nor do we require a
specific procedure to estimate them. We only impose high-level
conditions on the accuracy and complexity of the first step
estimate. In particular, our main result holds irrespectively
of whether the function $r_0$ is, for example, a density, a conditional mean
function or a quantile regression function, or whether it is
estimated by kernel methods, orthogonal series or sieves.
Moreover, our results are not confined to nonparametrically generated
covariates, but also apply in settings where $r_0$ is
estimated using parametric or semiparametric restrictions.

Our main result uses techniques from empirical process theory to show
that the presence of generated covariates affects the first-order
asymptotic properties of $\hat{m}_{LL}$ only through a
\emph{smoothed} version of the estimation error $\hat{r}(s)-r_0(s)$.
This additional smoothing typically improves the rate of convergence
of the estimator's stochastic part, reducing the ``curse of dimensionality''
from estimating $r_0$ to a secondary concern in this context.
It does not, however, affect the order of magnitude of the deterministic
component. Still, the estimator $\hat{m}_{LL}$ can have a faster
overall rate
of convergence than the first step estimator $\hat{r}$ if the latter has
a sufficiently small bias.

We extensively illustrate the implications of our main result for the important
special case that $r_0$ is the conditional mean function in an
auxiliary nonparametric regression. For this setting, we derive simple and
explicit stochastic expansions that can not only be used to
establish asymptotic normality or the rate of consistency of
the estimated regression function itself, but also study the
properties of more complex estimators, in which estimation of a regression
function merely constitutes an intermediate step, such as
structured nonparametric models imposing additive separability
[\citet{stone1985additive}]. Our results thus cover a wide range of models,
and should therefore be of general interest. We use our techniques
to study two such examples in greater detail: nonparametric estimation
of a
simultaneous equation model and nonparametric estimation of a
censored regression model.

To the best of our knowledge, there are only few papers on nonparametric
regression with estimated covariates not tailored to a specific application.
\citet{andrews1995nonparametric} derives some results for generated covariates
converging at a parametric rate. \citet{sperlich2009note} uses restrictive
assumptions which lead to asymptotic results that are different from the
ones obtained in the present paper. \citet{song2008uniform} considers series
estimation of the functional\vadjust{\goodbreak} $g(x,r) = \E(Y|r(X)=x)$ indexed by
$x\in\mathcal{X}\subset\IR$ and $r\in\Lambda$, where $\Lambda$
is a function space with finite integral bracketing entropy,
and derives a rate of consistency uniformly over $(x,r)\in\mathcal{X}
\times\Lambda$; see also \citet{einmahl2000uniform} for a~related
problem.

Our paper is also related to a recent literature on semiparametric
estimation problems with generated covariates. \citet{li2002semiparametric}
consider a partial linear model with generated covariates. \citet{hahn2010}
use pathwise derivatives to derive the influence function of semiparametric
linear GMM-type estimators. \citet{escanciano2011uniform} provide stochastic
expansions for sample means of weighted semiparametric regression residuals
with potentially generated regressors, and study their application to
certain index models.
Compared to the nonparametric problems studied in this paper, semiparametric
applications typically exhibit several additional technical issues. In
particular,
different techniques are needed to control the magnitude of certain
remainder terms.
Addressing these issues would require substantial refinements our
results, which
are not needed for the class of nonparametric problems we are focusing on.
To keep the present paper more readable, we study semiparametric
estimators with generated covariates separately in \citet{mammen2011semi}.

The outline of this paper is as follows. In the next section,
we describe our setup in detail. Section \ref{sec3} gives some motivating
examples. Section \ref{sec4} establishes the asymptotic theory and
states the main results. In Section~\ref{sec5}, we apply our results to
some of the examples given in Section \ref{sec3}, thus illustrating their
application in practice. Finally, Section \ref{sec6} concludes. All
proofs are collected in the \hyperref[app]{Appendix}.

\section{Nonparametric regression with generated covariates}

The nonparametric regression model with generated regressors can be
written as
\begin{equation}\label{model1}
Y = m_0(r_0(S)) + \varepsilon \qquad \mbox{with }\E( \varepsilon
|r_0(S))=0,
\end{equation}
where $Y$ is the dependent variable, $S$ is a $p$-dimensional
vector of covariates, $m_0\dvtx \IR^d\to\IR$ and $r_0\dvtx
\IR^p\to\IR^d$ are unknown functions and $\varepsilon$ is an
error term that has mean zero conditional on the true value of
covariates to covariates~$r_0(S)$.\footnote{Note that in contrast
to an earlier working paper version of this paper, we do no longer
assume that the ``index'' $r_0(S)$ is a sufficient statistic for
the covariates $S$, which would imply that $\E(Y|r_0(S))=\E(Y|S)$.}
We assume that there is additional information available
outside of the basic model~\eqref{model1} such that the
function~$r_0$ is identified. For example,~$r_0$ could be (some
known transformation of) the mean function in an auxiliary
nonparametric regression, which might involve another random
vector, say~$T$, in addition to~$Y$ and $S$.

Our aim is to estimate the function $m_0(x) = \E(Y|r_0(S)=x)$. Since
$r_0$ is unobserved,
obtaining a direct estimator based on a nonparametric regression of $Y$
on $R = r_0(S)$
is clearly not feasible. We therefore consider the following two-stage
procedure. In the
first stage, an estimate $\hat{r}$ of $r_0$ is obtained. We do not require
a specific estimator for this step. Instead, we only impose the
high-level restrictions that
the estimator $\hat{r}$ is uniformly consistent, converging at a rate
specified below, and takes
on values in a function class that is not too complex. Depending on the
nature of the
function $r_0$, these kind of regularity conditions are typically
satisfied by various common nonparametric
estimators, such as kernel-based procedures or series estimators, under
suitable smoothness restrictions.
In the second step, we then obtain our estimate $\hat m_{LL}$ of
$m_0$ through a nonparametric regression
of $Y$ on the generated covariates $\hat{R} = \hat{r}(S)$, using
local linear smoothing.
That is, our estimator is given by $\hat m_{LL}(x)=\hat\alpha
$ obtained from
\[
(\hat\alpha, \hat\beta) = \argmin_{\alpha,\beta}\sum
_{i=1}^n \bigl (Y_i- \alpha- \beta^T (\hat R_i - x)\bigr)^2 K_h(\hat R_i -x),
\]
where $ K_h(u) = \prod_{j=1}^d \mathcal{K}(u_j/h_j)/h_j$ is
a $d$-dimensional product kernel with univariate kernel function
$\mathcal{K}$,
and $h=(h_1,\ldots,h_d)$ is a vector of bandwidths that tend to zero as
the sample size $n$
increases to infinity.

For the later asymptotic analysis, it will also be useful to compare
$\hat m_{LL}$
to an infeasible estimator $\tilde m_{LL}$ that uses the true
function $r_0$
instead of an estimate $\hat{r}$. Such an estimator can be obtained by local
linear smoothing of~$Y$ versus $R=r_0(S)$, that is, it is given by
$\tilde m_{LL}(x)=\tilde\alpha$,
where
\[
(\tilde\alpha, \tilde\beta) = \argmin_{\alpha,\beta}\sum
_{i=1}^n \bigl(Y_i- \alpha- \beta^T (R_i - x)\bigr)^2 K_h(R_i -x).
\]
In order to distinguish these two estimators, we refer to $\hat
m_{LL}$ in the following
as the \emph{real} estimator, and to $\tilde m_{LL}$ as the \emph
{oracle} estimator.

Our use of local linear estimators in this paper is based on the
following considerations.
First, in a classical setting with fully observed covariates,
estimators based on local
linear regression are known to have attractive properties with regard
to boundary
bias and design adaptivity [see \citet{fan1996local} for an extensive
discussion], and
they allow a complete asymptotic description of their distributional properties.
In the present setting with generated covariates, these properties
simplify the asymptotic treatment.
The design adaptivity leads to a discussion of bias terms that does not
require regular densities for the randomly perturbed covariates, and
the complete asymptotic theory allows a clear description
of how the final estimator is affected by the estimation of the covariates.
On the other hand, our assumptions on the estimation of the covariates
are rather general and
can be verified for a broad class of smoothing methods, including
sieves and orthogonal series estimators.

\section{Motivating examples}\label{sec3}

There are many statistical applications which involve
nonparametric estimation of a regression \vadjust{\goodbreak}function using
nonparametrically generated covariates. In this section,
we give an overview of some of the most popular examples
and explain how they fit into our framework. In Section \ref{sec4},
we revisit the first three of these examples, studying their
asymptotic properties in detail. A thorough treatment of
the remaining examples involves several additional technical issues
beyond dealing with the presence of estimated covariates,
such as boundary problems, and is thus omitted for brevity.
See also \citet{mammen2011semi} for an extensive discussion
of semiparametric problems with generated covariates.

\subsection{The generic example: Nonparametric two-stage regression}

In many applications, the unknown function $r_0$ is a conditional
expectation function from an auxiliary nonparametric regression.
As a first motivating example, we therefore consider a ``two-stage''
nonparametric regression model given by
\begin{eqnarray}
 Y &=& m_0(r_0(S)) + \varepsilon\nonumber, \\
 T &=& r_0(S) + \zeta\nonumber, \nonumber
\end{eqnarray}
where $\zeta$ is an unobserved error term that satisfies
$E[\zeta| S] =E[\varepsilon| r_0(S)] = 0$. As the structure of
this example is particularly simple, it is used extensively in
Section~\ref{sec4} below to illustrate the application of our main result.
Proceeding like this is instructive, as the types of technical
difficulties encountered in this example are representative for
those in a wide range of other statistical applications.

\subsection{Nonparametric censored regression}
Consider a nonparametric regression model with fixed censoring, that is,
\begin{equation}\label{censor}
Y = \max\bigl(0,\mu_0(X) - U\bigr),
\end{equation}
where $U$ is an unobserved mean zero error term that is assumed
to be independent of the covariates $X$. Fixed censoring is a
common phenomenon in many applications, for example, the analysis of
wage data. Note that the censoring threshold could be different
from zero, as long as it is known.
\citet{linton2002nonparametric} establish identification of the
function $\mu_0$ under the tail condition
$\lim_{u\rightarrow-\infty} uF_U(u)=0$ on the distribution
function $F_U$ of~$U$. In particular, they show that the
function $\mu_0$ can be written as
\begin{equation}\label{censest}
\mu_0(x) = \lambda_0 - \int_{r_0(x)}^{\lambda_0} \frac{1}{q_0(r)}\,dr,
\end{equation}
where $r_0(x) = \E(Y|X=x)$, $q_0(r) = \E(\mathbb{I}\{Y>0\}|r_0(X) =r)$,
and $\lambda_0$ is some suitably chosen constant. An estimate
of the function $\mu_0$ can then be obtained from a sample
analog of \eqref{censest}, that is, through numerical integration
of a nonparametric estimate of the function $q_0(r)^{-1}$.
Nonparametric estimation of $q_0$ involves nonparametrically
generated regressors, and thus fits into our framework with
$(Y,S) = (\mathbb{I}\{Y>0\},X)$ and $r_0(S)=r_0(X)$.

\subsection{Nonparametric triangular simultaneous equation models}

Covariates that are correlated with disturbance terms appear in
many economic models and are denoted as endogenous. When, for example,
analyzing the relationship between wages and schooling,
unobserved individual characteristics like ability or
motivation might affect both the outcome and the explanatory
variable. A common approach is to model these quantities
jointly, achieving identification by using so-called
instrumental variables, that are independent of unobservables,
affect the endogenous variable, but exert no direct influence
on the outcome. Consider, for example, the nonparametric triangular
simultaneous equation model discussed in
\citet{newey1999nonparametric}, which is of the form
\begin{eqnarray}
\label{npsim1} Y &=& \mu_1(X_1,Z_1) + U ,\\
\label{npsim2} X_1 &=& \mu_2(Z_1,Z_2) + V .
\end{eqnarray}
Here the interest is in estimating the function $\mu_1$. To
achieve identification, one imposes the restrictions
$\E(V|Z_1,Z_2) = 0$, $\E(U) = 0$ and $\E(U|Z_1,Z_2,V) =
\E(U|V)$, which follow, for example, if the vector of exogenous
covariates and instruments $Z=(Z_1,Z_2)$ is jointly independent
of the disturbances $(U,V)$. Now let $m(x_1,z_1,v)
=\E(Y|X_1=x_1,Z_1=z_1,V=v)$. Under the above assumptions, it is
straightforward to show that
\[ \label{add}
m(x_1,z_1,v)=\mu_1(x_1,z_1) + \lambda(v),
\]
where $\lambda(v) = \E(U|V=v)$. The first component of this additive
model could, for example, be estimated by marginal integration
[\citet{newey1994kernel}, \citet{linton1995kernel}], which relies on the
fact that
\begin{equation}\label{eqmi}
\int m(x_1,z_1,v)f_V(v)\,dv = \mu_1(x_1,z_1),
\end{equation}
where $f_V$ is the probability density function of $V$. Implementing
a sample version of \eqref{eqmi} requires estimating the function $m$.
Since the residuals $V$ are not directly observed but must be estimated
by some
nonparametric method, this fits into our framework with $(Y,S) =
(Y,(X_1,Z_1,Z_2),X_1)$ and
$r_0(S) = (X_1,Z_1,X_1 - \mu_2(Z_1,Z_2))$.

\begin{remark}
An alternative to marginal integration would be an approach
based on smooth backfitting [\citet{mammen1999backfitting}].
Smooth backfitting estimators avoid several problems
encountered by mar\-ginal integration in case of covariates with
moderate or high dimension, but involves a more involved
statistical analysis which is beyond the scope of the present
paper. We are going to study smooth backfitting with nonparametrically
generated covariates in a separate paper.
\end{remark}

\subsection{Generalized Roy model}
D'Hautfoeuille and Maurel (\citeyear{dHautfoeuille2009Roy}) consider a generalized Roy
model of occupational choice that is related to the previous
example in the sense that it also leads to an additive\vadjust{\goodbreak}
regression model. Let $Y_k$ denote the individual's potential
earnings in sector $k\in\{0,1\}$ of an economy,
$X=(X_0,X_1,X_c)$ a vector of covariates, and assume that
$\E(Y_k|X,\eta_1,\eta_2)=\psi_k(X_k,X_c) + \eta_k$, where
$(\eta_0,\eta_1)$ are sector-specific productivity terms known
by the agent but unobserved by the analyst. Expected utility
from working in sector $k$ is assumed to be $U_k =
\E(Y_k|X,\eta_1,\eta_2) + G_k(X)$, the sum of sector-specific
expected earnings and a nonpecuniary component that depends on~$X$.
Along with $X$, the analyst observes the chosen sector
$D$, which satisfies $D=\mathbb{I}\{U_1 > U_0\},$ and the realized
earnings $Y=DY_1 + (1-D)Y_0$.

One object of interest in this context is the pair of functions $(\psi
_1, \psi_0)$.
Under some weak additional conditions, \citet{dHautfoeuille2009Roy}
show that
\[
\E(Y|D=d,X)= \psi_d(X_d,X_c) + \lambda_d\bigl(\Pr(D=d|X)\bigr)
\]
for $d\in\{0,1\}$, which is again an additive model involving
unobserved covariates,
namely the conditional probabilities $\Pr(D=d|X)$ of choosing sector~$d$. This setting
fits into our framework in the same way as the previous example.

\subsection{Nonparametric nonseparable triangular simultaneous
equation models}
Imbens and Newey (\citeyear{imbens2009identification}) consider a generalized version
of the above-mentioned triangular simultaneous equation model with nonadditive
disturbances:
\begin{eqnarray}
Y &=& \mu_1(X_1, Z_1, U),\\
X_1 &=& \mu_2(Z_1, Z_2, V).
\end{eqnarray}
Nonseparable models have become popular in the recent
econometric literature, as they allow for substantially more
general forms of unobserved heterogeneity than specifications
in which the disturbance terms enter additively. The focus here
is typically on averages of the function $\mu_1$, such as the
average structural function,
\[
\ASF(x_1,z_1) =
\E_U(\mu_1(x_1,z_1,U)).
\]
To achieve identification, assume
that the function $\mu_2$ is strictly monotone in its last
argument, that $V$ is continuously distributed, and that the
unobserved disturbances $(U,V)$ are jointly independent of $Z$.
Then it can be shown that $U$ and $(X_1,Z_1)$ are independently
conditional on the so-called control variable
$W=F_{X_1|Z}(X_1,Z)$, where $F_{X_1|Z}$ denotes the
distribution function of $X_1$ given $Z$. Under an additional
support condition, this result implies that the ASF is
identified through the relationship
\begin{equation}\label{asf}
\ASF(x_1,z_1) = \int m(x_1,z_1,w)\,dF_W,
\end{equation}
where $m(x_1,z_1,w) = \E(Y| X_1=x_1, Z_1=z_1, W=w)$.
Since the control variable $W$ is unobserved and has to be estimated in order
to implement a~sample analog\vadjust{\goodbreak} estimator of \eqref{asf}, this setting also
fits into the framework of this paper. In particular, nonparametric estimation
of $m$ is covered with $(Y,S) = (Y,(X_1,Z_1,Z_2),X_1)$ and
$r_0(S) = (X_1,Z_1,F_{X_1|Z}(X_1,Z))$.

\section{Asymptotic properties}\label{sec4}

It is straightforward to show that $\hat{m}_{LL}$
consistently estimates the function $m_0$ under standard
conditions. Obtaining refined asymptotic properties, however,
requires more involved arguments. In this section, we derive a
stochastic expansion of the difference between the real and the
oracle estimator, in which the leading terms are
kernel-weighted averages of the first stage estimation error.
This is our main result. It can be used, for example, to obtain uniform
rates of consistency for the real estimator, or to prove its
asymptotic normality. We demonstrate this in the
next section for specific forms of $r_0$ and~$\hat{r}$.

Throughout this section, we use the notation that for any
vector $a\in\IR^d$ the value $a_{\mathrm{min}}= \min_{1\leq j\leq d}
a_j$ denotes the smallest of its elements, $a_+ = \sum_{j=1}^d
a_j$ denotes the sum of its elements,
$a_{-k}=(a_1,\ldots,a_{k-1},a_{k+1},\ldots,a_d)$ denotes the
$d-1$-dimensional subvector of $a$ with the $k$th element
removed and $a^b = (a_1^{b_1},\ldots,a_d^{b_d})$ for any
vector $b\in\IR^d$. For ease of presentation in the following,
we avoid logarithmic terms in rates of convergence; that is, we
state assumptions and results in the form $o_P(n^\xi)$ instead
of $O_P(\log{n}^\gamma)$ with $\xi,\gamma>0$.

\subsection{Assumptions}

In order to analyze the asymptotic properties of the local
linear estimator with nonparametrically generated regressors,
we make the following assumptions.

\begin{assumption}[(Regularity conditions)]\label{as1} We assume the following properties
for the data distribution, the bandwidth, and kernel function $\mathcal{K}$:
\begin{longlist}[(iii)]
\item[(i)] The sample observations $(Y_i,S_i)$ are i.i.d.
\item[(ii)] The random vector $R=r_0(S)$ is continuously
distributed with compact support $I_R$. Its density
function $f_R$ is twice continuously differentiable
and bounded away from zero on $I_R$.
\item[(iii)] The function $m_0$ is twice continuously
differentiable on $I_R$.
\item[(iv)] $E[\exp(l | \varepsilon|) | S] \leq C$
almost surely for a constant $C>0$ and $l > 0$
small enough.
\item[(v)] The kernel function $\mathcal{K}$ is a twice
continuously differentiable, symmetric density
function with compact support, say $[-1,1]$.
%
\item[(vi)] The bandwidths $h=(h_1,\ldots,h_d)$ satisfies
$h_j \sim n^{-\eta_j}$ for $j=1,\ldots,d$ and
$\eta_+ < 1$.
\end{longlist}
\end{assumption}

Most conditions in Assumption \ref{as1} are standard regularity and smoothness
conditions for kernel-type nonparametric regression, with the exception
of Assumption~\ref{as1}(iv). The subexponential tails of $\varepsilon$ conditional
on $S$ assumed there are needed to apply certain results from empirical
process theory in our proofs. Such a condition is not very restrictive though.\vadjust{\goodbreak}

\begin{assumption}[(Accuracy)]\label{as2} The components $\hat r_j$ and $ r_{0,j}$
of $\hat r$ and $ r_{0}$, respectively, satisfy
\[
\sup_{s}| \hat r_j(s) - r_{0,j}(s)| = o_P(n^{-\delta_j})
\]
for some $\delta_j > \eta_j$ and all $j=1,\ldots,d$.\vspace*{-2pt}
\end{assumption}

Assumption \ref{as2} is a ``high-level'' restriction on the accuracy of the
estimator $\hat{r}$.
It requires each component of the estimate of the function $r_0$ to be uniformly
consistent, converging at rate at least as fast as the corresponding
bandwidth in the second
stage of the estimation procedure. This is typically not a restrictive
condition, and it
allows for estimators $\hat{r}$ that converge at a~rate slower than
the oracle estimator
$\tilde m_{LL}$. Uniform rates of consistency are widely available
for all common
nonparametric estimators; see, for example, \citet
{masry1996multivariate} for results
on the Nadaraya--Watson, local linear and local polynomial estimators, or
\citet{newey1997convergence} for series estimators.\vspace*{-2pt}

\begin{assumption}[(Complexity)]\label{as3} There exist sequences of sets ${\mathcal
M}_{n,j}$ such
that:
\begin{longlist}[(ii)]
\item[(i)]$\Pr(\hat r_j \in{\mathcal M}_{n,j})\to1$ as $n\to\infty
$ for all $j=1,\ldots,d$.
\item[(ii)] For a constant $C_M>0$ and a function $r_{n,j}$ with $\|r_{n,j} -
r_{0,j}\|_\infty=o( n^{-\delta_j})$, the set ${{\overline{\mathcal{
M}}}}_{n,j}={\mathcal M}_{n,j}\cap\{r_j\dvtx \|r_j - r_{n,j}\|_\infty
\leq n^{-\delta_j}\}$ can be covered by at most $C_M \exp(\lambda
^{-\alpha_j} n ^{\xi_j})$ balls with $\|\cdot\|_{\infty}$-radius
$ \lambda$ for all $ \lambda\leq n^{-\delta_j}$, where $0<\alpha_j
\leq2$, $\xi_j \in\IR$ and $\|\cdot\|_{\infty}$ denotes the
supremum norm.\vspace*{-2pt}
\end{longlist}
\end{assumption}

Assumption \ref{as3} requires the first-stage estimator $\hat{r}$ to
take values in a function space ${\mathcal M}_{n,j}$ that is not
too complex, with probability approaching 1. Here the
complexity of the function space is measured by the cardinality
of the covering sets. This is a typical requirement for many
results from empirical process theory; see \citet{van1996weak}.
The second part of Assumption \ref{as3} is
typically fulfilled under suitable smoothness restrictions. For
example, suppose that ${\mathcal M}_{n,j}$ is the set of functions
defined on some compact set $I_S \subset\mathbb{R}^p$ whose
partial derivatives up to order $k$ exist and are uniformly
bounded by some multiple of $n^{\xi_j^*}$ for some
$\xi_j^*\geq0$. Then Assumption \ref{as3}(ii) holds with $\alpha_j =
p/k$ and $\xi_j=\xi_j^*\alpha_j$ [\citet{van1996weak}, Corollary
2.7.2]. For kernel-based estimators of $r_0$, one
can then verify part (i) of Assumption \ref{as3} by explicitly
calculating the derivatives. Consider, for example, the one-dimensional
Nadaraya--Watson estimator $\hat r_{n,j}$ with bandwidth of
order $n^{-1/5}$. Choose $r_{n,j}$ equal to $ r_{0,j}$ plus
asymptotic bias term. Then one can check that the second
derivative of $\hat r_{n,j}- r_{n,j}$ is absolutely
bounded by $O_P(\sqrt{\log n}) = o_P(n^{\xi_j^*})$ for all
$\xi_j^*>0$. For sieve and orthogonal\vspace*{1pt} series estimators,
Assumption \ref{as3}(i) immediately holds when the set ${\mathcal M}_{n,j}$
is chosen as the sieve set or as a subset of the linear span of
an increasing number of basis functions, respectively.
For a discussion of entropy bounds and further references, we
refer to \citet{vandegeer2009book}.\vadjust{\goodbreak}

\begin{assumption}[(Continuity)]\label{as4}
For any $r \in\mathcal{M}_n=\mathcal{M}_{n,1}\times\cdots
\times\mathcal{M}_{n,d}$ the conditional expectation $\tau^B(x,r) =
\E(\rho(S)|
r(S)= x )$ with $\rho(S)=\E(Y|S) - \E(Y|r_0(S))$ exists and is twice
differentiable with respect to its first argument, with derivatives
that are
uniformly bounded in absolute value, and satisfies
\[
\|\tau^B(x,r_1)-\tau^B(x, r_2)\| \leq C^*_B \|r_1-r_2\|_{\infty
}\qquad\mbox{a.s.}
\]
for all $r_1,r_2\in\mathcal{M}_n$
and a constant $C^*_B>0$.
\end{assumption}

Assumption \ref{as4} imposes certain smoothness restrictions on the conditional
expectation of $\rho(S)$. The term $\rho(S)$ can be thought of as capturing
the influence of the underlying covariates $S$ on the outcome variable $Y$
that is not excreted through the ``index'' $r_0(S)$. In certain applications,
the ``index''~$r_0(S)$ is a sufficient statistic for the function
$m_0$, and thus
$\rho(S)=0$ with probability 1. In this case, Assumption \ref{as4} is
trivially satisfied.
Note that $\rho(S)=\E(\varepsilon|S)$, and that $\tau^B(\cdot
,r_0)\equiv0$ by construction.

\subsection{The key stochastic expansion}

With the assumptions given in the previous section, we are now
ready to state our main result, which is a~stochastic expansion of the real
estimator $\hat m_{LL}(x)$ around the oracle estimator $\tilde
m_{LL}(x)$.
Our aim is to derive an explicit characterization of the influence of
the presence
of generated regressors on the final estimator of the function $m_0$.
To this end, we define $ w(x,r) =(1, (r_1(S) - x_1)/h_1,\ldots,(r_d(S)
- x_d)/h_d)$,
and set $N_h(x) = \E( w(x,r) w(x,r)^T K_h( r(S)-x))$. Next, we define
\begin{eqnarray*}
\Delta(x,r) &=& e_1^{\top}N_h(x)^{-1}\E\bigl( K_h\bigl(r_0(S)-x\bigr)w(x,r) \bigl(r(S) -r_0(S)\bigr)\bigr),\\
\Gamma(x, r) &=& e_1^{\top}N_h(x)^{-1}\E\bigl( K_h'\bigl(r_0(S)-x\bigr)^\top w(x,r)\bigl(r(S) -r_0(S)\bigr)\rho(S)\bigr)
\end{eqnarray*}
for any $r \in\mathcal{M}_n$, where $K^{\prime}_h(u) =(
\mathcal{K}_{h,j}^{\prime}(u)\dvtx j=1,\ldots,d)^\mathsf{T}$ is a
vector with elements
$\mathcal{K}_{h,j}^{\prime}(u)=\mathcal{K^{\prime}}(u_{j}/h_{j})/h^2_{j}
\prod_{j^*\not=j}\mathcal{K}(u_{j^*}/h_{j^*})/h_{j^*}$. Finally, we put
$\hat\Delta(x) = \Delta(x,\hat r)$ and $\hat\Gamma(x) = \Gamma
(x,\hat r)$.
With this notation, we can now state our main theorem.

\begin{theorem}\label{mainexpansion}Suppose Assumptions \ref{as1}--\ref{as4} hold. Then
\begin{eqnarray*}
\sup_{x \in \mathrm{I_R}} | \hat m_{LL}(x) - \tilde m_{LL}(x)
+m_0^\prime(x) \hat\Delta(x) -\hat\Gamma(x) |= O_P(n ^{-\kappa}),
\end{eqnarray*}
where
$\kappa= \min\{\kappa_1,\ldots,\kappa_3\}$
with
\begin{eqnarray*}
\kappa_1 &<& {\frac{1}{2}} (1 - \eta_+) + (\delta-\eta)_{\mathrm{min}} -
\frac{1}{2} \max_{1\leq j \leq d} (\delta_j
\alpha_j + \xi_j), \\
\kappa_2 &<& 2 \eta_{\mathrm{min}} + (\delta-\eta)_{\mathrm{min}},\\
\kappa_3 &<& \delta_{\mathrm{min}}+(\delta-\eta)_{\mathrm{min}}.
\end{eqnarray*}
\end{theorem}

The two leading terms in our stochastic expansion of the real estimator~$\hat m_{LL}(x)$ around
the oracle estimator $\tilde m_{LL}(x)$, which are accounting for
the presence of generated covariates,
are both smoothed versions of the first-stage estimation error
$\hat r(s) - r_0(s)$. To see this more
clearly, note that it follows from standard arguments for local
polynomial smoothing that
\begin{eqnarray*}
\Delta(x,r) &=& \frac{ \E( K_h(r_0(S) - x) (r(S) - r_0(S)))}{f_R(x)}
+ O_P(n ^{-\kappa})\qquad \mbox{and}\\
\Gamma(x,r) &=& \frac{ \E( K_h'(r_0(S)-x)^\top(r(S) - r_0(S))\rho
(S_i))}{f_R(x)} + O_P(n ^{-\kappa}),
\end{eqnarray*}
uniformly over $x\in I_{R,n}^-=\{x\in I_R \dvtx$ the support of $K_h(\cdot
-x)$ is a subset of~$I_R\}$.
In order to achieve a certain rate of convergence for the real
estimator, it is thus not necessary to
have an estimator of $r_0$ that converges with the same rate or a
faster one, since the asymptotic
properties of the estimator using nonparametrically generated
regressors only depend on a~smoothed version of the first-stage estimation error. While smoothing
does not affect the order
of the estimator's deterministic part, it typically reduces the
variance and thus allows for less precise
first-stage estimators. Note that the first adjustment term is
negligible in regions where the
regression function is flat, since $m_0'(x)=0$ in this case.
Conversely, the impact of generated
covariates is accentuated when the true regression function is steep.
Also note that $\hat\Gamma(x)=0$
when $\E(\varepsilon|S)=0$, as the latter implies that $\rho
(s)\equiv0$. This is a~natural condition
in certain empirical applications.\vspace*{-2pt}

\begin{remark}
In Theorem \ref{mainexpansion} no assumptions are made about the process
generating the data for estimation of $r_0$. In particular,
nothing is assumed about dependencies between the errors in the
pilot estimation and the regression errors $\varepsilon_i$. We
conjecture that better rates than $n^{-\kappa}$ can be proven
under such additional assumptions, but the results would only
be specific to the respective full model under consideration.
One way to extend our approach to such a setting would be to
use our empirical process methods to bound the remainder term
of higher order differences between $\hat m$ and
$\tilde m$, and to treat the leading terms of the resulting higher
order expansion by other, more direct methods.\vspace*{-2pt}
\end{remark}

\section{Examples revisited}\label{sec5}

In this section, we apply our high-level results from Section \ref{sec4}
to some of the motivating examples presented in Section \ref{sec3},
which are representative for the others in terms of employed
techniques. Assuming a specific nature of the function $r_0$
and a specific method to estimate it,\vspace*{1pt} explicit
forms of the adjustment terms $\hat\Delta(x)$ and $\hat\Gamma(x)$
in Theorem \ref{mainexpansion} can be
derived in order to account for the presence of generated
covariates. Our focus in this section is on the practically most
important case that $r_0$ is the conditional mean function in an auxiliary
nonparametric regression. Many other applications can be
treated along the same lines.\vadjust{\goodbreak}

\subsection{Generic example: Two-stage nonparametric regression}
The main setting in which we illustrate the application of the stochastic
expansion from Theorem \ref{mainexpansion} is the ``two-stage'' nonparametric regression model
given by\looseness=-1
\begin{eqnarray*}
 Y&=&m_0(r_0(S)) + \varepsilon ,\\
 T&=&r_0(S) + \zeta,
\end{eqnarray*}\looseness=0
where $\zeta$ is an unobserved error term that satisfies
$E[\zeta| S] =E[\varepsilon| r_0(S)] = 0$. For simplicity, we focus
on the case that $R=r_0(S)$ is
a one-dimensional covariate, but generalizations to multiple generated
covariates or the presence
of additional observed covariates are immediate.

Our strategy for deriving asymptotic properties of $\hat m_{LL}$
in this framework is to first provide an explicit representation for the
adjustment terms $\hat\Delta(x)$ and $\hat\Gamma(x)$ from Theorem \ref{mainexpansion},
which are then combined with standard results about the oracle
estimator $\tilde m_{LL}$. For this approach it is convenient to use
a~kernel-based smoother to estimate $r_0$. Since the bias of both $\hat
\Delta(x)$ and
$\hat\Gamma(x)$ is of the same order as of this first-stage estimator,
we propose to estimate the function $r_0$ via $q$th order local
polynomial smoothing, which includes the local linear estimator as
the special case $q=1$. Formally, the estimator is given by $\hat
{r}(s) =
\hat\alpha$, where
\begin{equation}\qquad
(\hat\alpha,\hat\beta) = \argmin_{\alpha,\beta} \sum_{i=1}^n \biggl(
T_i - \alpha- \sum_{1 \leq u_+ \leq q}\beta_{r}^T (S_{i}-s)^u\biggr)^2 L_g(S_i-s)
\end{equation}
and $L_g(s) = \prod_{j=1}^p \mathcal{L}(s_j/g)/g$ is a
$p$-dimensional product kernel\vspace*{1pt} built from the univariate kernel
$\mathcal{L}$, $g$ is a vector of bandwidths, whose components
are assumed to be the same for simplicity, and $\sum_{1 \leq
u_+ \leq q}$ denotes the summation over all
$u=(u_1,\ldots,u_p)$ with $1\leq u_+ \leq q$. When $r_0$ is
sufficiently smooth,\vspace*{1pt} the asymptotic bias of local polynomial
estimators of order $q$ is well known to be $O(g^{q+1})$
uniformly over $x\in I_R$ (if $q$ is uneven), and
can thus be controlled. A further technical advantage of using
local polynomials is that the corresponding estimator admits a
certain stochastic expansion under general conditions, which is
useful for our proofs. We make the following assumption, which
is essentially analogous to Assumption \ref{as1}, except for Assumption~\ref{as4}(iii). This additional assumption requires higher order
smoothness of the kernel, necessary to bound the $k$th
derivative of the estimator $\hat r$. This allows us to verify
Complexity Assumption \ref{as3} for $\hat r$.

\begin{assumption}\label{as5} We assume the following properties
for the data distribution, the bandwidth and kernel function $\mathcal{L}$:
\begin{longlist}[(iii)]
\item[(i)] The observations $(S_i,Y_i,T_i)$ are i.i.d., and the
random vector $S$ is continuously distributed with
compact support $I_S$. Its density function
$f_S$ is bounded and bounded away from zero on $I_S$. It
is also differentiable with a bounded derivative. The residuals $\zeta$
satisfy $\E|\zeta|^{\epsilon} < \infty$ for some $\epsilon>2$.\vadjust{\goodbreak}
\item[(ii)] The function $r_0$ is $q+1$ times continuously
differentiable on $I_S$.
\item[(iii)] The kernel function $\mathcal{L}$ is a $k$-times
continuously differentiable, symmetric density
function with compact support, say $[-1,1]$, for some
natural number $k\geq\max\{2, p/2\}$.
\item[(iv)]The bandwidth satisfies $g \sim n^{-\theta}$ for
some $0<\theta<1/p$.\vspace*{-2pt}
\end{longlist}
\end{assumption}

To simplify the presentation, we also assume that the function
$r_0(s)$ is strictly monotone in at least one of its arguments,
which can be taken to be the last one without loss of
generality. This assumption could be easily removed at the cost
of a substantially more involved notation in the following
results.\vspace*{-2pt}

\begin{assumption}\label{as6} The function $r_0(u_{-p},u_p)$ is strictly monotone
in $u_p$, and we have that $r_0(u_{-p}, \varphi(u_{-p},x))=x$ for
some twice continuously differentiable function $\varphi$.\vspace*{-2pt}
\end{assumption}

The following proposition shows that in the present context
the function $\hat\Delta(x)$ can be written as the sum of a
smoothed version of the first stage estimator's bias function,
a kernel-weighted average of the first-stage residuals
$\zeta_1,\ldots,\zeta_n$, and some higher order remainder
terms. For a concise presentation of the result, we introduce
some particular kernel functions. Let $L^*$ denote the
$p$-dimensional equivalent kernel of the local polynomial
regression estimator, given in \eqref{equivkernel} in the
\hyperref[app]{Appendix}, and define the one-dimensional kernel functions
\begin{eqnarray*}
J_h(x,s) &=& \int K_h\bigl(r_0(s)-x - \partial_s r_0(s) uh\bigr) L^*(u)\, du, \\
H^\Delta_{g}(x,v) &=& \frac{\partial_{x}\varphi(v_{-p},x)}{g} \int L^*\biggl(s_{-p},\frac{\varphi(v_{-p},x)-v_p}{g} + s_p \partial
_{-p}\varphi(v_{-p},x)\biggr)\, ds.\vspace*{-2pt}
\end{eqnarray*}
%
Then, with this notation, we obtain the following proposition.

\begin{proposition} \label{proptwostage} Suppose that Assumptions
\ref{as1} and \ref{as4}--\ref{as6} hold. Then we have for the correction factor $\hat\Delta$ in
Theorem \ref{mainexpansion} that
\[
\sup_{x\in I_{R}}|\hat\Delta(x) - \hat\Delta_A(x) - \hat\Delta
_B(x)| = O_p\biggl(\frac{\log(n)}{ng^p} \biggr),
\]
where the terms $\hat\Delta_A(x)$ and $\hat\Delta_B(x)$ satisfy
\begin{eqnarray*}
\sup_{x\in I_{R}}|\hat\Delta_A(x)| &=& O_p\bigl(\bigl(\log(n)/(n\max\{g,h\}
)\bigr)^{1/2}\bigr) \quad \mbox{and} \\
 \sup_{x\in I_{R}}|\hat\Delta_B(x)| &=& O_p(g^{q+1}).
\end{eqnarray*}
Moreover, uniformly over
$x\in I_{R,n}^-$, it is $\hat\Delta_B(x) = g^{q+1} E[ b(S)|
r_0(S)=x]+ o_p(g^{q+1})$ with a bounded function $b(s)$ given
in \eqref{fsbias} in the \hyperref[app]{Appendix}, and the term\vadjust{\goodbreak} $\hat
\Delta_A(x)$ allows for the following expansions uniformly over
$x\in I_{R,n}^-$, depending on the limit of $g/h$:
\begin{longlist}[(a)]
\item[(a)] If $g/h\to0$, then
\[
\hat\Delta_A(x) 
= \frac{1}{n f_R(x)}\sum_{i=1}^n K_h\bigl(r_0(S_i)-x\bigr) \zeta_i +
O_p\biggl(\biggl(\frac{g^2}{h^2}+ \frac{g^{3/2}}{h}\biggr)\biggl(\frac{\log(n)}{nh}\biggr)^{1/2}\biggr).
\]
%
%
\item[(b)] If $h=g$, then
\[
\hat\Delta_A(x)
= \frac{1}{nf_R(x)} \sum_{i=1}^n J_h(x, S_i) \zeta_i +O_p\biggl(\biggl(\frac
{\log(n)}{n}\biggr)^{1/2}\biggr).
\]
\item[(c)] If $g/h\to\infty$, then
\[
\hat\Delta_A(x) 
= \frac{1}{nf_R(x)} \sum_{i=1}^n H^{\Delta}_g(x, S_i) \zeta_i +
O_p\biggl(\frac{g^2}{h^2}\biggl(\frac{\log(n)}{ng}\biggr)^{1/2}+\biggl(\frac{\log(n)}{n}\biggr)^{1/2}\biggr).
\]
\end{longlist}
\end{proposition}

It should be emphasized that in all three cases of the above
proposition the leading term in the expression for $\hat
\Delta_A(x)$ is equal to an average of the error terms
$\zeta_i$ weighted by a \emph{one-dimensional} kernel function,
irrespective of $p=\dim(S)$. The dimension of the covariates
thus affects the properties of $\hat\Delta(x)$ only through
higher-order terms. Furthermore, it should be noted that one
can also derive expressions of $\hat\Delta(x)$ similar to the
ones above for values of $x$ close to the boundary of the
support. Likewise these take the form of a one-dimensional
kernel weighted average of the error terms $\zeta_i$ plus a
higher-order term. The corresponding kernel function, however,
has a more complicated closed form varying with the point of
evaluation.

The following proposition establishes a result similar to
Proposition \ref{proptwostage} for the second adjustment term $\hat\Gamma(x)$.
We again introduce a particular one-dimensional kernel function,
defined as
\[
H_g^{\Gamma}(x,v)
%
=\int g^{-1} L^*\biggl(s_{-p},\frac{\varphi(v_{-p},x)-v_p}{g} + s_p
\,\partial_{p}\varphi(v_{-p},x)\biggr)\,ds
\lambda(v_{-p},x)
\]
with
\[
\lambda(v_{-p},x) =\frac{\partial_{v_p} (\rho(v_{-p},\varphi
(v_{-p},x) f_S(v_{-p},\varphi(v_{-p},x))
\det(\partial_{v_{-p}}\varphi(v_{-p},x))}{f_S(v_{-p},\varphi
(v_{-p},x))\,\partial_{v_p}r_{0}(v_{-p},\varphi(v_{-p},x))},
\]
where $L^*$ still denotes the
$p$-dimensional equivalent kernel of the local polynomial
regression estimator, given in \eqref{equivkernel} in the
\hyperref[app]{Appendix}.

\begin{proposition} \label{proptwostage2b} Suppose that Assumptions
\ref{as1}
and \ref{as4}--\ref{as6} hold. Then we have that
\[
\sup_{x\in I_{R}}|\hat\Gamma(x) - \hat\Gamma_A(x) - \hat\Gamma
_B(x)| = O_p\biggl(\frac{\log(n)}{ng^p} \biggr),
\]
where the terms $\hat\Gamma_A(x)$ and $\hat\Gamma_B(x)$ satisfy
\[
\sup_{x\in I_{R}}|\hat\Gamma_A(x)| = O_p\bigl(\bigl(\log(n)/(ng)\bigr)^{1/2}\bigr)
\quad \mbox{and}\quad\sup_{x\in I_{R}}|\hat\Gamma_B(x)| = O_p(g^{q+1}).
\]
Moreover, uniformly over $x\in I_{R,n}^-$, it is $\hat\Gamma_B(x) =
g^{q+1}\, \partial_x E[ b(S)\rho(S)|
r_0(S)=x]+ o_p(g^{q+1})$ with a bounded function $b(s)$ given
in \eqref{fsbias} in the \hyperref[app]{Appendix}, and the term $\hat
\Gamma_A(x)$ allows for the following expansion uniformly over
$x\in I_{R,n}^-$:
\begin{equation}
\hat\Gamma(x) = \frac{1}{n f_R(x)}\sum_{i=1}^n H_g^{\Gamma
}(x,S_i)\zeta_i+ o_P\Biggl(\sqrt{\frac{\log(n)}{ng}}\Biggr).
\end{equation}
\end{proposition}

Again, the leading term in the expression for $\hat\Gamma_A(x)$ is
equal to an
average of the error terms $\zeta_i$ weighted by a \emph
{one-dimensional} kernel
function, and thus behaves similarly to one-dimensional nonparametric regression
estimator. A~similar result could be established for regions close to
the boundary
of the support. Note that in contrast to Proposition \ref{proptwostage}, the details
of the result in Proposition \ref{proptwostage2b} do not depend
on the relative magnitude of the bandwidths used in the first and
second stage
of the estimation procedure.

Combining Theorem \ref{mainexpansion} and Propositions \ref{proptwostage}--\ref{proptwostage2b}
with well-known results about the oracle estimator $\tilde
{m}_{LL}$, various asymptotic properties
of the real estimator $\hat{m}_{LL}$ can be derived. In the
following corollaries we present results
for the most relevant scenarios, addressing uniform rates of
consistency and stochastic
expansions of order $o_P(n^{-2/5})$ for proving pointwise asymptotic normality.
More refined expansions of higher orders such as $o_P(n^{-1/2})$, which are
useful for the analysis of semiparametric problems in which $m_0$ plays
the role of an infinite
dimensional nuisance parameter [e.g., \citet{newey1994variance}, \citet{andrews1994asymptotics}, \citet{chen2003estimation}],
would also be possible. We do not present such results here as they
would require strong smoothness restrictions that are unattractive in
applications.
See \citet{mammen2011semi} for an alternative approach to controlling
the influence of generated
covariates in semiparametric models.

Starting with considering the uniform rate of consistency, it is
well known [\citet{masry1996multivariate}] that under
Assumption \ref{as1} the oracle estimator satisfies
\[
\sup_{x\in I_R}|\tilde{m}_{LL}(x) - m(x)| = O_p\bigl(\bigl(\log(n)/nh\bigr)^{1/2}
+h^2\bigr).
\]
This implies the following result.

\begin{corollary}\label{corAtwostage}Suppose that Assumptions \ref{as1}, \ref{as4} and
\ref{as5} hold.
Then
\[
\sup_{x\in I_R}|\hat{m}_{LL}(x) - m(x)| = O_p\biggl(\frac{\log
(n)^{1/2}}{ (n\max\{h,g\})^{1/2}} + h^2 +\frac{\log
(n)}{ng^p}+g^{q+1} + n^{-\kappa}\biggr).
\]
\end{corollary}
%

Straightforward calculations show that, under appropriate smoothness
restrictions, it is possible to recover the oracle rate for the real
estimator given
suitable
choice of $\eta$ and $\theta$, even if the first-stage estimator
converges at a~strictly slower rate. Note that the rate in Corollary \ref{corAtwostage} improves
upon a~bound on the uniform rate of convergence of a two-stage
regression estimator
derived in \citet{ahn1995nonparametric} for a similar setting.

Next, we derive stochastic expansions of $\hat{m}_{LL}$ of
order $o_P(n^{-2/5})$ for the case that $\eta=1/5$. Such expansions
immediately imply results on pointwise asymptotic normality of the
real estimator. We start with the case that $\theta= \eta$, in which
the stochastic terms $\hat\Gamma_A(x)$ and $\hat\Delta_A(x)$
are of the same order of magnitude (other bandwidth choices will be discussed
below). During the analysis of this setting, it becomes clear
that applying Theorem \ref{mainexpansion} requires $p \theta< 3/10$.
Thus in order to use the expansion in Proposition \ref{proptwostage}(b), only $p=1$ is
admissible;
that is, $S$ must be one-dimensional for the choice $\theta= \eta$
to be feasible. In this setting, the notation for the kernel functions
appearing in the stochastic expansions can be somewhat simplified.
We define
\begin{eqnarray*}
\tilde{J}(v,x)&=&\int K\bigl(v-r_0^{\prime}(r_0^{-1}(x))u\bigr) L^*(u)\,du,\\[-4pt]
\tilde H^\Gamma(v,x) &=& \int L^*\bigl(v + s\, \partial_{x}r_0^{-1}(x)\bigr)\,ds\tilde\lambda(x),
\end{eqnarray*}
where
\[
\tilde\lambda(x) =\frac{\partial_{v} (\rho(r_0^{-1}(x))
f_S(r_0^{-1}(x)))}{f_S(r_0^{-1}(x))r_0^{\prime}(r_0^{-1}(x))},
\]
where $r_0^{-1}$ is the inverse function of $r_0$, which exists by
Assumption \ref{as6}.\vspace*{-4pt}

\begin{corollary} \label{corBtwostage} Suppose that Assumptions \ref{as1} and
\ref{as4}--\ref{as6} hold with $\eta=\theta=1/5$ and $p=q=1$.
Then the following expansions hold uniformly over $x\in I_{R,n}^-$:
\begin{eqnarray*}
&&\hat{m}_{LL}(x) -m_0(x)\\[-3pt]
&&\qquad =\frac{1}{nf_R(x)} \sum_{i=1}^n
K_h\bigl(r_0(S_i)-x\bigr) \varepsilon_i \\[-3pt]
 && \qquad \quad {} - \frac{1}{nf_R(x)} \sum_{i=1}^n \bigl(m_0'(x) \tilde
{J}_h\bigl(r_0(S_i)-x,x\bigr) - \tilde H_h^\Gamma\bigl(S_i-r_0^{-1}(x),x\bigr)\bigr) \zeta_i \\[-3pt]
 && \qquad \quad {} + \frac{1}{2} \beta(x)h^2 + o_p(n^{-2/5}) ,
\end{eqnarray*}
where the bias is given by
\begin{eqnarray*}
\beta(x) &=& \int u^2 K(u)\,du m_0^{\prime\prime}(x) \\[-3pt]
&& {} - \int u^2 L(u)\,du\bigl(r_0^{\prime\prime} (r_0^{-1}(x))m_0'(x) - \partial
_x[r_0^{\prime\prime} (r_0^{-1}(x))\rho(r_0^{-1}(x))] \bigr).\vadjust{\goodbreak}
\end{eqnarray*}
In particular, we have
\[
(nh)^{1/2}\bigl(\hat{m}_{LL}(x) - m_0(x)- \beta(x)h^2\bigr) \stackrel
{d}{\rightarrow} N(0,\sigma_m^2(x)),
\]
where
$\sigma_m^2(x) = [ \Var(\varepsilon|R=x)\int K(t)^2\,dt - 2
E(\varepsilon\zeta|R=x)\int K(t)(\tilde J(t,x) m_0'(x) -
\tilde{H}^\Gamma(t,x))\,dt \Var(\zeta|R=x)\int(m_0'(x) \tilde
J(t,x)- \tilde{H}^\Gamma(t,x))^2\,dt]/f_R(x) $
is the asymptotic variance.
\end{corollary}

Under the conditions of the corollary, the limiting distribution
of $\hat{m}_{LL}(x)$ is generally affected by the pilot estimation
step, although a
qualitative description of the impact seems difficult. Depending on the
curvature
of $m_0$ and the covariance of $\varepsilon$ and $\zeta$, the
asymptotic variance
of the estimator using generated regressors can be bigger or smaller
than that of
the oracle estimator $\tilde{m}_{LL}$. There thus exist settings
where in practice it would be preferable to base inference on the real
estimator even if
one was actually able to compute the oracle estimator.

The next corollary considers the case that $\theta>\eta$, and thus
$g/h \to0$.
Again, applying Theorem \ref{mainexpansion} requires $p \theta< 3/10$ in this setting, and
thus only $p=1$ is admissible when using Proposition \ref{proptwostage}(a) for such a
choice of
bandwidths. The corollary also focuses on the special case that $\rho
(S) :=
\E(Y|R) - \E(Y|S)=0$, which implies that $\hat\Gamma(x) = 0$
with probability 1. This condition is satisfied for certain empirical
applications,
such as, for example, models IV models. Without this additional
restriction, an expansion
of the difference $\hat{m}_{LL}(x)-m_0(x)$ would be dominated by
the term $\hat\Gamma_A(x)$,
which is $O_p((\log(n)/(ng))^{1/2})$ and thus converges at a \emph
{slower} rate than
the oracle estimator.

\begin{corollary} Suppose that Assumptions \ref{as1}, \ref{as4} and \ref{as5} hold with $\eta
=1/5$, $1/5< \theta< 3/10$ and $p=q=1$,
and that $\rho(S)=0$ with probability 1.
Then the following expansion holds uniformly over $x\in I_{R,n}^-$:
\begin{eqnarray*}
\hat{m}_{LL}(x)-m_0(x) &=& \frac{1}{n f_R(x)}\sum_{i=1}^n
K_h\bigl(r_0(S_i)-x\bigr) \bigl(\varepsilon_i - m_0^{\prime}(x)\zeta_i\bigr) \\
&& {} + {\frac{1}{2}} h^2 \int u^2 K(u)\,du m_0^{\prime\prime} (x) +
o_p(n^{-2/5}) .
\end{eqnarray*}
In particular, we have
\[
(nh)^{1/2}\biggl(\hat{m}_{LL}(x) - m_0(x)- \frac{1}{2} h^2 \int u^2
K(u)\,du m_0^{\prime\prime} (x)\biggr) \stackrel{d}{\rightarrow}
N(0,\sigma_m^2(x)),
\]
where
$\sigma_m^2(x) = \Var(\varepsilon- m_0'(R)\zeta|R=x)\int K(t)^2\,dt /
f_R(x) $
is the asymptotic variance.
\end{corollary}

The limiting distribution of $\hat{m}_{LL}(x)$ is again affected by
the use of generated covariates under the conditions of the corollary.
In this particular case, the\vadjust{\goodbreak} form of the asymptotic variance has an
intuitive interpretation: the estimator $\hat{m}_{LL}(x)$ has
the same limiting distribution as the local linear oracle estimator in the
hypothetical regression model
\[
Y = m_0(r_0(S)) + \varepsilon^*,
\]
where $\varepsilon^* = \varepsilon- m_0'(r_0(S))\zeta$.
As in Corollary \ref{corBtwostage} above, depending on the curvature of $m_0$ and the
covariance of $\varepsilon$ and $\zeta$, the asymptotic variance of the
estimator using generated regressors can be bigger or smaller than that
of the oracle estimator $\tilde{m}_{LL}$.

The next corollary discusses the case when $\theta< \eta$.
For such a choice of bandwidth, applying Theorem \ref{mainexpansion} requires
no restrictions on the dimensionality of $S$. It turns out
that in this case $\hat{m}_{LL}(x)= \tilde{m}_{LL}(x) +
o_p(n^{-2/5})$,
and thus the limit distribution of $\hat{m}_{LL}$ is
the same as for the oracle estimator~$\tilde{m}_{LL}$. The
effect exerted by the presence of nonparametrically generated
regressors is thus first-order asymptotically negligible for conducting
inference on $m_0$ in this case.

\begin{corollary} \label{corCtwostage} Suppose that Assumptions \ref{as1},
\ref{as4} and \ref{as5} hold with
$\theta< \eta=1/5$. Then the following expansion
holds uniformly over $x\in I_{R,n}^-$ if $\frac{2}{5}(q+1)^{-1} <
\theta< \frac{3}{10} p^{-1} $:
\begin{eqnarray*}
\hat{m}_{LL}(x) - m_0(x) &=& \frac{1}{n f_R(x)}\sum_{i=1}^n
K_h\bigl(r_0(S_i)-x\bigr) \varepsilon_i \\
&& {} + \frac{1}{2}h^2 \int u^2 K(u)\,du
m_0^{\prime\prime} (x) + o_p(n^{-2/5}) .
\end{eqnarray*}
In particular, we have
\[
(nh)^{1/2}\biggl(\hat{m}_{LL}(x) - m_0(x)- \frac{1}{2} h^2 \int u^2
K(u)\,du m_0^{\prime\prime} (x)\biggr) \stackrel{d}{\rightarrow}
N(0,\sigma_m^2(x)),
\]
where
$\sigma_m^2(x) = \Var(\varepsilon|R=x)\int K(t)^2\,dt / f_R(x) $
is the asymptotic variance.
\end{corollary}

\subsection{Nonparametric censored regression}
Consider estimation of the censored regression model
in \eqref{censor}. Let $\hat{r}(x)$ be the $q$th order local
polynomial estimator of the conditional mean $r_0(x)=
\E(Y|X=x)$, and let $\hat{q}(r)$ be the local linear estimator
of $q_0(r)$ using the generated covariates $\hat{r}(X_i)$. Then
an estimate of $\mu_0$ is given by
%
%
\begin{equation}\label{censestimator}
\hat{\mu}(x)= \lambda+\int_{\hat{r}(x)}^{\lambda} \frac
{1}{\hat{q}(u)}\,du,
\end{equation}
where the constant $\lambda$ is chosen large enough to satisfy
$\lambda>\max_{i=1,\ldots,n}\hat{r}(X_i)$ with probability
tending to one. Generalizing \citet{linton2002nonparametric},
we consider the use of higher-order local polynomials for the
first stage estimator, and allow the bandwidth used for the
computation of $\hat{r}$ and $\hat{q}$ to be different. For
presenting the asymptotic properties of $\hat{\mu}$, let
$s_0(x)=\E(\mathbb{I}\{Y>0\}|X=x)$ be the proportion of uncensored
observations conditional\vadjust{\goodbreak} on $X=x$, and assume that this
function is continuously differentiable and bounded away from
zero on the support of $X$. We then obtain the following
result.

\begin{corollary}\label{corcens}
Suppose that Assumptions \ref{as1} and \ref{as5} hold with $(Y,S,T) = (\mathbb{I}\{
Y>0\},X,Y)$ and
$R=r_0(S)=r_0(X)$. Furthermore, suppose that
$\theta\in(\underline{\theta} ,\bar{\theta})$ where
$\underline{\theta}$ and $\bar{\theta}$ are constants depending
on $\eta$, $q$ and $p$ as follows:
\[
\bar{\theta} = \frac{1-3\eta}{p}
\quad\mbox{and}\quad
\underline{\theta}= \max\biggl\{\frac{1-4\eta}{p},\frac{1}{2(q+1)+p}\biggr \}.
\]
Under these conditions, we have that
\[
\sqrt{ng^p}\bigl(\hat{\mu}(x) - \mu_0(x)\bigr) \stackrel{d}{\rightarrow}
N\biggl(0,\frac{\sigma^2_r(x)}{f_S(x)s_0^2(x)}\int L(t)^2\,dt\biggr),
\]
where $\sigma^2_r(x)= \Var(Y|X=x)$.
\end{corollary}

The corollary is analogous to Theorem 5 in
\citet{linton2002nonparametric}. However, using our results,
substantially simplifies the proof and provides insights on
admissible choices of bandwidths. Note that the lower bound
$\underline{\theta}$ is chosen such that both the bias of
$\hat{r}$ and $\hat{q}$ tends to zero at a rate faster than
$(ng^p)^{-1/2}$. Due to this undersmoothing, the limiting
distribution of $\hat{\mu}-\mu$ is centered at zero. Note that
the final estimator converges at the same rate as the generated
regressors. This is due to the fact that the function $\hat{r}$
is not only used to compute $\hat{q}$, but also determines the
limits of integration in~\eqref{censestimator}. The ``direct''
influence of the generated regressors in the estimation of $q$
is asymptotically negligible in this particular application.

\subsection{Nonparametric triangular simultaneous equation models}

Now consider nonparametric estimation of the structural
function $\mu_1$ in the triangular simultaneous equation
model \eqref{npsim1}--\eqref{npsim2} using a marginal integration
estimator. In order to keep the notation simple, we
restrict our attention to the arguably most relevant case with
a single endogenous regressor, but allow for an arbitrary
number of exogenous regressors and instruments. Let~$\hat{\mu}_2(z)$ be the $q$th order local polynomial estimator
of $\mu_2(z) = \E(X_1|Z=z)$, and let $\hat{m}(x_1,z_1,v)$ be
the local linear estimator of $m(x_1,z_1,v) =\E(Y|X_1
=x_1,Z_1=z_1,V=v)$. The latter is\vspace*{2pt} computed using the generated
covariates $\hat{V}_i = X_{1i}-\hat{\mu}_2(Z_i)$ instead of the
true residuals $V_i$ from
equation \eqref{npsim2}. For simplicity, we use the same bandwidth for
all components of $\hat{m}$; that is, we
put $\eta_j \equiv\eta$ for all $j=1,\ldots,(2+d_1)$. The marginal
integration
estimator of $\mu_1(x_1,z_1)$ is then given by the following sample
version of \eqref{eqmi}:
\begin{equation}\label{marint}
\hat\mu_1(x_1,z_1) = \frac{1}{n}\sum_{i=1}^n \hat{m}(x_1,z_1,\hat
{V}_i).
\end{equation}
The following result establishes the estimator's asymptotic
normality.

\begin{corollary}\label{corsim}
Suppose that Assumption \ref{as1} holds with $(Y,S,T) =
(Y,(X_1,\allowbreak Z_1,Z_2),X_1)$ and $R=r_0(S) = (X_1,Z_1, X_1
-\mu_2(Z_1,Z_2))$, and that Assumption~\ref{as5} holds with $r_0(S) =
\mu_2(Z_1,Z_2)$. Furthermore, suppose that $\eta\in
( \max\{1/\allowbreak(5+d_1),1/(2p+3)\}, 1/(1+d_1))$, and that
$\theta\in(\underline{\theta},\bar{\theta})$, where
$\underline{\theta}$ and $\bar{\theta}$ are constants depending
on $\eta$, $q$ and $d_j=\dim(Z_j)$ as follows:
\[
\bar{\theta} = \frac{1-3\eta}{2p}\quad\mbox{and} \quad\underline{\theta}=
\frac{1-\eta(d_1+1)}{2(q+1)} ,
\]
where $p=d_1+d_2$. Under these conditions, we have that
\[
\sqrt{nh^{1+d_1}}\bigl(\hat{\mu}_1(x_1,z_1) - \mu_1(x_1,z_1) \bigr) \stackrel
{d}{\rightarrow} N\biggl(0, \E\biggl(\frac{\sigma_\varepsilon
^2(x_1,z_1,V)}{f_{XZ|V}(x_1,z_1,V)}\biggr)
\int\tilde{K}(t)^2\,dt\biggr),
\]
where $\tilde{K}(t) = \prod_{i=1}^{1+d_1}\mathcal{K}(t_i)$ is a
$(1+d_1)$-dimensional product kernel, and
$\sigma_\varepsilon^2(x_1, z_1, v) = \Var(Y- m(R)|R=(x_1,z_1,v))$.
\end{corollary}

Under the conditions of the corollary, the asymptotic
variance of $\hat\mu_1(x_1,z_1)$ is not influenced by the
presence of generated regressors: If $\hat{m}$ was replaced
in \eqref{marint} with an oracle estimator $\tilde{m}$
using the actual disturbances $V_i$ instead of the
reconstructed ones, the result would not change. Also,
note that the exclusion restrictions on the instruments
imply that $\E(Y|X_1,Z_1,V) = \E(Y|X_1,Z_1,Z_2)$.\vspace*{1pt} Therefore
Assumption \ref{as4} is automatically satisfied, and the adjustment
term $\hat\Gamma(x)$ from Theorem~\ref{mainexpansion} is equal to zero and does not
have to be considered for the proof.

\section{Conclusions}\label{sec6}

In this paper, we analyze the properties of nonparametric
estimators of a regression function, when some the covariates
are not directly observable, but have been estimated by a
nonparametric first-stage procedure. We derive a stochastic
expansion showing that the presence of generated regressors
affects the limit behavior of the estimator only through a
smoothed version of the first-stage estimation error. We apply
our results to a number of practically relevant statistical
applications.

\begin{appendix}\label{app}
\section*{Appendix: Proofs}

Throughout the Appendix, $C$ and $c$ denote generic constants
chosen sufficiently large or sufficiently small, respectively,
which may have different values at each appearance.
Furthermore, define $\bar{\mathcal{M}}_n =
\bar{\mathcal{M}}_{n,1} \times\cdots\times
\bar{\mathcal{M}}_{n,d}$.

\subsection{\texorpdfstring{Proof of Theorem \protect\ref{mainexpansion}}{Proof of Theorem 1}}
In order to prove the statement of the theorem,
we have to introduce some notation. Throughout the proof of this and
the following
statements, we denote the unit vector $(1,0,\ldots,0)^{T}$ in
$\mathbb{R}^{p+1}$ by $e_1$. We also write $ w_i(x,r) =(1, (r_1(S_i) -
x_1)/h_1,\ldots,(r_d(S_i) - x_d)/h_d)$,
and put $ w_i(x)=w_i(x,r_0)$, $ \hat w_i(x)=w_i(x,\hat r)$ and $ \tilde
w_i(x)=w_i(x,\tilde r)$.
We also\vspace*{1pt} define $M_h(x,r) = n^{-1}\sum_{i=1}^n w_i(x,r) w_i(x,r)^T K_h(
r(S_i)-x)$,
and put $M_h(x)=M_h(x, r_0)$, $\hat M_h(x)=M_h(x,\hat r)$ and\vadjust{\goodbreak}
$\tilde M_h(x)=M_h(x,\tilde r)$
and set $N_h(x) = \E( M_h(x,\allowbreak r_0))$.
Furthermore, define $\varepsilon^*= \varepsilon- \rho(S)$ and note
that we have
$\E(\varepsilon^*|S)=0$ by construction. It also holds
that
\[
Y_i = m_0(r_0(S_i)) + \varepsilon_i^* + \rho(S_i).
\]
Next, it follows from standard calculations that the real estimator
$\hat m_{LL}$ can be written as
\[
\hat m_{LL}(x) = m_0(x) + \hat m_{LL,A}(x)+ \hat
m_{LL,B}(x)+ \hat m_{LL,C}(x)+ \hat m_{LL,D}(x)+ \hat m_{LL,E}(x),
\]
where $\hat m_{LL,j}(x) = \hat\alpha_j$ for
$j\in\{A,B,C,D,E\}$, and
\begin{eqnarray*}
(\hat\alpha_A, \hat\beta_A) &=& \argmin_{\alpha,\beta
}\sum_{i=1}^n \bigl(\varepsilon_i^*- \alpha- \beta^T \bigl(\hat r(S_i) -
x\bigr)\bigr)^2 K_h\bigl(\hat r(S_i) -x\bigr),\\[-2pt]
(\hat\alpha_B, \hat\beta_B) &=& \argmin_{\alpha,\beta
}\sum_{i=1}^n \bigl(m_0(r_0(S_i)) - m_0(x)- m_0^{\prime}(x)^T \bigl(r_0(S_i) -
x\bigr) \\ &&\hspace*{144pt}{} -\alpha- \beta^T \bigl(\hat r(S_i) - x\bigr)\bigr)^2\\[-2pt]
&&\hphantom{\argmin_{\alpha,\beta
}\sum_{i=1}^n}{}\times K_h\bigl(\hat r(S_i)-x\bigr),\\[-2pt]
(\hat\alpha_C, \hat\beta_C) &=& \argmin_{\alpha,\beta
}\sum_{i=1}^n \bigl(-m_0^{\prime}(x)^T\bigl(\hat{r}(S_i)-r_0(S_i)\bigr) - \alpha-
\beta^T \bigl(\hat r(S_i) - x\bigr)\bigr)^2 \\[-2pt]
&&\hphantom{\argmin_{\alpha,\beta
}\sum_{i=1}^n}{}\times K_h\bigl(\hat r(S_i) -x\bigr), \\[-2pt]
(\hat\alpha_D, \hat\beta_D) &=& \argmin_{\alpha,\beta
}\sum_{i=1}^n \bigl(m_0^{\prime}(x)^T\bigl(\hat{r}(S_i)-x\bigr)- \alpha- \beta
^T \bigl(\hat r(S_i) - x\bigr)\bigr)^2 \\[-2pt]
&&\hphantom{\argmin_{\alpha,\beta
}\sum_{i=1}^n}{}\times K_h\bigl(\hat r(S_i) -x\bigr), \\[-2pt]
(\hat\alpha_E, \hat\beta_E) &=& \argmin_{\alpha,\beta
}\sum_{i=1}^n \bigl(\rho(S_i)- \alpha- \beta^\mathsf{T} \bigl(\hat
r(S_i) - x\bigr)\bigr)^2 K_h\bigl( \hat r(S_i) - x\bigr).
\end{eqnarray*}
Similarly, the oracle estimator $\tilde m_{LL}$ can be
represented as
\[
\tilde m_{LL}(x) = m_0(x) + \tilde m_{LL,A}(x)+ \tilde
m_{LL,B}(x)+\tilde m_{LL,C}(x) +\tilde m_{LL,D}(x) + \tilde
m_{LL,E}(x),
\]
where $\tilde m_{LL,j}(x) = \tilde\alpha_j$ for
$j\in\{A,B,C,D,E\}$, and
\begin{eqnarray*}
(\tilde\alpha_A, \tilde\beta_A) &=& \argmin_{\alpha,\beta
}\sum_{i=1}^n \bigl(\varepsilon_i- \alpha- \beta^T \bigl(r_0(S_i) - x\bigr)\bigr)^2
K_h\bigl(r_0(S_i) -x\bigr),\\[-2pt]
(\tilde\alpha_B, \tilde\beta_B) &=& \argmin_{\alpha,\beta
}\sum_{i=1}^n \bigl(m_0(r_0(S_i)) - m_0(x)- m_0^{\prime}(x)^T \bigl(r_0(S_i) -x\bigr) \\[-2pt]
 && \hspace*{140pt} {} - \alpha- \beta^T \bigl(r_0(S_i) - x\bigr)\bigr)^2 \\[-2pt]
 &&\hphantom{\argmin_{\alpha,\beta
}\sum_{i=1}^n}{}\times K_h\bigl(r_0(S_i) -x\bigr),\\[-2pt]
(\tilde\alpha_C, \tilde\beta_C) &=& \argmin_{\alpha,\beta
}\sum_{i=1}^n \bigl(-m_0^{\prime}(x)^T\bigl(\hat r(S_i) - r_0(S_i)\bigr)- \alpha
- \beta^T \bigl(r_0(S_i) - x\bigr)\bigr)^2 \\[-2pt]
 &&\hphantom{\argmin_{\alpha,\beta
}\sum_{i=1}^n}{}\times K_h\bigl(r_0(S_i) -x\bigr) \\[-2pt]
(\tilde\alpha_D, \tilde\beta_D) &=& \argmin_{\alpha,\beta
}\sum_{i=1}^n \bigl(m_0^{\prime}(x)^T \bigl(r_0(S_i) - x\bigr)- \alpha- \beta^T
\bigl(r_0(S_i) - x\bigr)\bigr)^2 \\[-2pt]
 &&\hphantom{\argmin_{\alpha,\beta
}\sum_{i=1}^n}{}\times K_h\bigl(r_0(S_i) -x\bigr). \\[-2pt]
(\tilde\alpha_E, \tilde\beta_E) &=& \argmin_{\alpha,\beta
}\sum_{i=1}^n \bigl(\rho(S_i)- \alpha- \beta^\mathsf{T} \bigl( r(S_i) -
x\bigr)\bigr)^2 K_h\bigl( r(S_i) - x\bigr).
\end{eqnarray*}
Note that by construction,
\begin{equation}\label{lemma01}
\hat m_{LL,D}(x) \equiv\tilde m_{LL,D}(x) \equiv0.
\end{equation}
We now argue that
\begin{equation}\label{lemma12}
\sup_{x\in I_R}|\hat m_{LL,A}(x) - \tilde m_{LL,A}(x)| =
O_p(n^{-\kappa_1}).
\end{equation}
For a proof of (\ref{lemma12}) note that $\hat m_{LL,A}(x)$
and $\tilde m_{LL,A}(x)$ are given by the first  elements of
the vectors $\hat M(x)^{-1} n^{-1} \sum_{i=1}^n
K_h(\hat r(S_i) - x) \varepsilon_i \hat w_i(x)$ and
$M(x)^{-1}\* n^{-1} \sum_{i=1}^n K_h( r_0(S_i) - x)
\varepsilon_i \tilde w_i(x)$, respectively. Using these
representations, one sees that (\ref{lemma12}) follows from Lemmas
\ref{lemAtheo1} and \ref{lemBtheo1} below.

As a second step, we now show that
\begin{equation}\label{lemma5}\qquad
\sup_{x\in I_R}|\hat m_{LL,E}(x) - \tilde m_{LL,E}(x) - \hat
\Gamma(x)| =O_p(n^{-\kappa_1}+n^{-\kappa_2}+n^{-\kappa_3}).
\end{equation}
To prove \eqref{lemma5}, put $\hat\mu(x) =\frac{1}{n} \sum
_{i=1}^n K_h(
\hat r(S_i) -x) \hat w_i(x) \rho(S_i)$ and $\mu(x)=\frac
{1}{n} \sum_{i=1}^n
K_h( r_0(S_i) -x)w_i(x) \rho(S_i)$, and write $G(x)= e_1^\mathsf{T}
(N_{h}(x))^{-1} \E(\hat\mu(x) - \mu(x))$. With this notation,
$\hat m_{LL,E}(x) = e_1^\mathsf{T}\hat M_{h}(x)^{-1} \hat
\mu(x)$ and
$ \tilde m_{LL,E}(x) = e_1^\mathsf{T}\times  M_{h}(x)^{-1} \mu(x)$. Using
Lemma \ref{lemEtheo1} and some results of Lemma \ref{lemCtheo1}, we
then find that
\begin{eqnarray*}
&&\hat m_{LL,E}(x) - \tilde m_{LL,E}(x) - G(x)\\[-2pt]
&& \qquad\!\! = e_1^\mathsf{T} \bigl(\hat M_{h}(x)^{-1} \hat\mu(x) -
M_{h}(x)^{-1} \mu(x) - \E( M_{h}(x))^{-1}\E\bigl(\hat\mu(x) - \mu
(x )\bigr)\bigr)\\[-2pt]
&& \qquad\!\! = O_P\bigl(n^{-((1/2)(1-\eta_{+})+(\delta-\eta)_{\mathrm{min}})}+
n^{-((1/2)(1-\eta_{+})+\delta_{\mathrm{min}})} +n^{-\kappa
_1}\bigr)=O_P(n^{-\kappa_1})
\end{eqnarray*}
uniformly over $x\in I_R$. Using standard smoothing arguments, we also
get that
\begin{eqnarray*}
G(x)&=& e_1^\mathsf{T} N_h(x)^{-1} \E\bigl(\hat\mu(x) - \mu(x)\bigr)\\[-2pt]
&=& \frac{1}{f_R(x)}\int\bigl(K_h\bigl( \hat r(u) -x\bigr) - K_h\bigl( r_0(u) -x\bigr)
\bigr)\rho(u) f_S(u)\,dx\,du\\[-2pt]
&& {} + O_P\bigl(n^{-2\eta_{\mathrm{min}}-(\delta-\eta)_{\mathrm{min}}}\bigr)\\[-2pt]
&=& \frac{1}{f_R(x)}\int K'_h\bigl( r_0(u) -x\bigr) \bigl( \hat r(u)- r_0(u)\bigr)
\rho(u) f_S(u)\,dx\,du\\[-2pt]
&& {} + O_P\bigl(n^{-\delta_{\mathrm{min}}-(\delta-\eta
)_{\mathrm{min}}}\bigr)+O_P(n^{-\kappa_2})\\[-2pt]
&=& \hat\Gamma(x)+ O_P(n^{-\kappa_2})+O_P(n^{-\kappa_3})
\end{eqnarray*}
uniformly over $x\in I_R$. This shows the claim in \eqref{lemma5}.

Finally, from Lemmas \ref{lemBtheo1} and \ref{lemCtheo1} we get that
\begin{eqnarray}
\label{lemma3}\sup_{x\in I_R}|\hat m_{LL,B}(x) - \tilde m_{LL,B}(x)| &=&
O_p(n^{-\kappa_2}),\\
\label{lemma4}\sup_{x\in I_R}|\hat m_{LL,C}(x) - \tilde m_{LL,C}(x)| &=&
O_p(n^{-\kappa_3}),
\end{eqnarray}
and it is easy to see that
\begin{equation}\label{lemma6}
\sup_{x\in I_R}|\tilde m_{LL,C}(x)-m_0'(x)\hat\Delta(x)| =
O_p(n^{-\kappa}).
\end{equation}
Taken together, the results in \eqref{lemma01}--\eqref{lemma6} imply
the statement of the theorem.

\begin{lemma} \label{lemAtheo1}
Suppose that the conditions of Theorem \ref{mainexpansion} hold. Then
\begin{eqnarray*}
&&\sup_{x\in I_R, r_1,r_2\in\bar{\mathcal{M}}_n} \Biggl| \frac{1}{n}\sum
_{i=1}^n K_h\bigl(r_1(S_i)-x\bigr)\varepsilon_i -\frac{1}{n}\sum_{i=1}^n
K_h\bigl(r_2(S_i)-x\bigr)\varepsilon_i\Biggr| \\
&& \qquad = O_p(n^{-\kappa_1}),\\
&&\sup_{x\in I_R, r_1,r_2\in\bar{\mathcal{M}}_n} \Biggl| \frac{1}{n}\sum
_{i=1}^n K_h\bigl(r_1(S_i)-x\bigr)\frac{r_{1,j}(S_i) -x_j}{h_j}\varepsilon_i \\
&&\hphantom{\sup_{x\in I_R, r_1,r_2\in\bar{\mathcal{M}}_n} \Biggl|}{}-\frac{1}{n}\sum_{i=1}^n K_h\bigl(r_2(S_i)-x\bigr)\frac{r_{2,j}(S_i)
-x_j}{h_j}\varepsilon_i\Biggr|\\
&& \qquad = O_p(n^{-\kappa_1}).
\end{eqnarray*}
\end{lemma}

\begin{pf}
We only prove the first statement of the lemma. The second
claim can be shown using essentially the same arguments.
Without loss of generality, we also assume that
\begin{equation}\label{k1}
\kappa_1 > (\delta- \eta)_{\mathrm{min}}.
\end{equation}
If $\kappa_1 \leq(\delta- \eta)_{\mathrm{min}}$ the statement of the
lemma follows from a direct bound. For $C_1,C_2>0$ large enough
(see below) we choose $C_{\varepsilon}$ such that
\begin{eqnarray}
\label{t5}\Pr\Bigl(\max_i|\varepsilon_i|> C_{\varepsilon} \log(n)\Bigr) &\leq&
n^{-C_1}, \\
\label{t6}\bigl|\E\varepsilon_i\mathbb{I}\{|\varepsilon_i| \leq C_{\varepsilon}
\log(n)\}\bigr|
&\leq& n^{-C_2}.
\end{eqnarray}
With this choice of $C_{\varepsilon}$ we define
\[
\Delta_i(r_1,r_2) = \bigl(K_h\bigl(r_1(S_i)-x\bigr) - K_h\bigl(r_2(S_i)-x\bigr) \bigr)\varepsilon_i^*
\]
with
\[
\varepsilon_i^* = \varepsilon_i \mathbb{I}\{|\varepsilon_i| \leq
C_{\varepsilon_i} \log(n)\} - \E\bigl(\varepsilon_i
\mathbb{I}\{|\varepsilon_i| \leq C \log(n)\}\bigr).
\]
For the proof of the lemma we apply a chaining argument; compare, for
example, the proof of Theorem 9.1 in
\citet{vandegeer2009book}. Now for $s\geq0$, let
$\bar{\mathcal{M}}_{s,n,j}^*$ be a set of functions chosen such
that for each $r\in\bar{\mathcal{M}}_{n,j}$ there exists
$r^*\in\bar{\mathcal{M}}_{s,n,j}^*$ such that $\|r-r^*\|_\infty
\leq2^{-s} n^{-\delta_j}$. That is, the functions in~$\bar{\mathcal{M}}_{s,n,j}^*$ are the midpoints of a $(2^{-s}
n^{-\delta_j})$-covering of $\bar{\mathcal{M}}_{n,j}$. By
Assumption~\ref{as3}, the set $\bar{\mathcal{M}}_{s,n,j}^*$ can be
chosen such that its cardinality $\#
\bar{\mathcal{M}}_{s,n,j}^*$ is at most $C
\exp((2^{-s}n^{-\delta_j})^{-\alpha_j} n ^{\xi_j})$.
Furthermore, define $\bar{\mathcal{M}}_{s,n}^* =
\bar{\mathcal{M}}_{s,n,1}^* \times\cdots\times
\bar{\mathcal{M}}_{s,n,d}^*$.

For $r_1,r_2 \in\bar{\mathcal{M}}_{n}$ we now choose
$r_1^s,r_2^s\in\bar{\mathcal{M}}_{s,n}^*$ such that $
\|r_{1,j}^s - r_{1,j}\|_{\infty} \leq2^{-s}n^{-\delta_j}$ and
$\|r_{2,j}^s - r_{2,j}\|_{\infty} \leq C 2^{-s}n^{-\delta_j}$,
for all $j$. We then consider the chain
\begin{eqnarray*}
\Delta_i(r_1,r_2) &=& \Delta_i(r_1^0,r_2^0) - \sum_{s=1}^{G_n}\Delta
_i(r_1^{s-1},r_1^s)
+ \sum_{s=1}^{G_n}\Delta_i(r_2^{s-1},r_2^s)\\[-2pt]
&&{}- \Delta
_i(r_1^{G_n},r_1) + \Delta_i(r_2^{G_n},r_2),
\end{eqnarray*}
where $G_n$ is the smallest integer that satisfies $G_n>
(1+c_G)(\kappa_1 - (\delta-\eta)_{\mathrm{min}})\*\log(n)/\log(2)$ for a
constant $c_G > 0$. With this choice of $G_n$, we obtain that
for $l=1,2$
\begin{equation} \label{t1}\qquad
T_1 = \Biggl|\frac{1}{n}\sum_{i=1}^n \Delta_i(r_l^{G_n},r_l) \Biggr| \leq C \log
(n) 2^{-G_n} n^{-(\delta-\eta)_{\mathrm{min}}} \leq C n^{-\kappa_1}.
\end{equation}
Now for any $a>c_G$ define the constant $c_a =
(\sum_{s=1}^\infty2^{-as})^{-1}$. It then follows that
\begin{eqnarray*}
 &&\Pr\Biggl(\sup_{r_1 \in\bar{\mathcal{M}}_n }\Biggl|\frac{1}{n}\sum_{i=1}^n
\sum_{s=1}^{G_n} \Delta_i(r_1^{s-1},r_1^s)\Biggr|>n^{-\kappa_1} \Biggr)\\[-2pt]
 && \qquad \leq\sum_{s=1}^{G_n} \Pr\Biggl(\sup_{r_1 \in\bar{\mathcal{M}}_n }
\Biggl|\frac{1}{n}\sum_{i=1}^n \Delta_i(r_1^{s-1},r_1^s) \Biggr|>c_a
2^{-as}n^{-\kappa_1} \Biggr) \\[-2pt]
&& \qquad \leq\sum_{s=1}^{G_n} \# \bar{\mathcal{M}}^*_{s-1,n} \# \mathcal
{\bar{M}}^*_{s,n}
\Pr\Biggl(\frac{1}{n}\sum_{i=1}^n \Delta_i(r_1^{*,s},r_1^{**,s}) >c_a
2^{-as}n^{-\kappa_1} \Biggr) \\[-2pt]
&& \qquad \quad {} + \sum_{s=1}^{G_n} \# \bar{\mathcal{M}}^*_{s-1,n} \# \mathcal
{\bar{M}}^*_{s,n}
\Pr\Biggl(\frac{1}{n}\sum_{i=1}^n \Delta_i(\tilde{r}_1^{*,s},\tilde
{r}_1^{**,s}) <c_a 2^{-as}n^{-\kappa_1} \Biggr) \\[-2pt]
&&\qquad  = T_2 + T_3,
\end{eqnarray*}
where the functions $r_1^{*,s},\tilde{r}_1^{*,s} \in
\bar{\mathcal{M}}^*_{s-1,n}$ and
$r_1^{**,s},\tilde{r}_1^{**,s}\in\bar{\mathcal{M}}^*_{s,n}$ are
chosen such that
\begin{eqnarray*}
&&\Pr\Biggl(\frac{1}{n}\sum_{i=1}^n \Delta_i(r_1^{*,s},r_1^{**,s}) >c_a
2^{-as}n^{-\kappa_1}\Biggr ) \\
&& \qquad = \max_{r_1^{s-1}, r_1^{s}} \Pr\Biggl(\frac{1}{n}\sum_{i=1}^n \Delta
_i(r_1^{s-1},r_1^{s}) >c_a 2^{-as}n^{-\kappa_1} \Biggr),\\
&&\Pr\Biggl(\frac{1}{n}\sum_{i=1}^n \Delta_i(\tilde{r}_1^{*,s},\tilde
{r}_1^{**,s}) <c_a 2^{-as}n^{-\kappa_1} \Biggr) \\
&& \qquad = \max_{r_1^{s-1}, r_1^{s}} \Pr\Biggl(\frac{1}{n}\sum_{i=1}^n \Delta
_i(r_1^{s-1},r_1^{s}) >c_a 2^{-as}n^{-\kappa_1} \Biggr).
\end{eqnarray*}
We now show that both $T_2$ and $T_3$ tend to zero at an
exponential rate:
\begin{eqnarray}
\label{t2} T_2 &\leq&\exp(-cn^c), \\
\label{t3} T_3 &\leq&\exp(-cn^c) .
\end{eqnarray}
We only show \eqref{t2}, as the statement \eqref{t3} follows by
essentially the same arguments. Using Assumption \ref{as3}, we obtain
by application of the Markov inequality that
\begin{eqnarray} \label{s2}
\qquad
T_2 &\leq& C \sum_{s=1}^{G_n} \prod_j \exp\bigl((2^{-s}n^{-\delta
_j})^{-\alpha_j}n^{\xi_j} \bigr)\nonumber\\
\qquad&&\hspace*{25pt}{}\times\E\Biggl(\exp\Biggl(\gamma_{n,s} \frac{1}{n}\sum_{i=1}^n \Delta
_i(r_1^{*,s},r_1^{**,s}) - \gamma_{n,s}c_a 2^{-as}n^{-\kappa_1} \Biggr)\Biggr)
\nonumber\\[-8pt]
\\[-8pt]
\qquad&\leq& C \sum_{s=1}^{G_n} \exp\biggl(\sum_j 2^{s\alpha_j} n^{\delta
_j\alpha_j +\xi_j} - \gamma_{n,s}c_a2^{-as}n^{-\kappa_1}\biggr)\nonumber\\
\qquad&&\hspace*{25pt}{}\times\prod_{i=1}^n \E\biggl( \exp\biggl(\gamma_{n,s} \frac{1}{n}\Delta
_i(r_1^{*,s},r_1^{**,s}) \biggr)\biggr),\nonumber
\end{eqnarray}
where $\gamma_{n,s} = c_\gamma2^{(2-a)s}n^{-\kappa_1+1-\eta_+
+ 2(\delta-\eta)_{\mathrm{min}}}$ with a constant $c_\gamma>0$, small
enough. Now the last term on the right-hand side of \eqref{s2}
can be bounded as follows:
\begin{eqnarray}\label{s1}\qquad\qquad
\E\biggl( \exp\biggl(\gamma_{n,s} \frac{1}{n}\Delta_i(r_1^{*,s},r_1^{**,s}) \biggr)\biggr)
&\leq&1 + C \E(\gamma_{n,s}^2n^{-2}\Delta
_i^2(r_1^{*,s},r_1^{**,s}))\nonumber
\\[-8pt]
\\[-8pt]
\qquad\qquad&\leq&\exp\bigl( C \gamma_{n,s}^2n^{-2} n^{\eta_+ - 2(\delta-\eta
)_{\mathrm{min}}}2^{-2s} \bigr),\nonumber
\end{eqnarray}
where we have used that
\begin{eqnarray*}
\biggl|\gamma_{n,s} \frac{1}{n}\Delta_i(r_1^{*,s},r_1^{**,s})\biggr|&\leq& C
\gamma_{n,s} \frac{1}{n}\log(n) n^{\eta_+}n^{-(\delta-\eta)_{\mathrm{min}}}2^{-s}\\
&\leq& C \log(n) n^{(\delta-\eta)_{\mathrm{min}}-\kappa_1} 2^{-as+s}\\
&\leq& C \log(n) n^{(c_G-a)(\kappa_1- (\delta-\eta)_{\mathrm{min}}) }\\
&\leq& C
\end{eqnarray*}
for $n$ large enough because of \eqref{k1}.
Inserting \eqref{s1} into \eqref{s2}, we obtain,
if $a$ and $c_\gamma$ were chosen sufficiently small, that
\begin{eqnarray*}
T_2  &\leq& C \sum_{s=1}^{G_n} \exp\biggl(\sum_j 2^{s\alpha_j} n^{\delta
_j\alpha_j +\xi_j} - c2^{2(1-a)s}n^{1-2\kappa_1-\eta_+ +2(\delta
-\eta)_{\mathrm{min}}} \biggr)\\
 &\leq& C \sum_{s=1}^{G_n} \exp( - c^{s}n^c )\\
 &\leq& \exp(-cn^c).
\end{eqnarray*}
Finally, it follows from a simple argument that
\begin{equation}\label{t4}\qquad
T_4 = \Pr\Biggl(\sup_{r_1,r_2 \in\bar{\mathcal{M}}_n }\Biggl|\frac{1}{n}\sum
_{i=1}^n \Delta_i(r_1^{0},r_2^0)\Biggl|>n^{-\kappa_1} \Biggr)\leq\exp(-cn^c)
\end{equation}
because the set $\bar{\mathcal{M}}_{0,n}^*$ can always be
chosen such that it contains only a~single element.

From \eqref{t1}, \eqref{t2}, \eqref{t3} and \eqref{t4}, we thus
obtain that
\begin{eqnarray}\label{t7}
\qquad&&\hspace*{22pt}\sup_{x\in I_R}\Pr\Biggl(\sup_{r_1,r_2 \in\bar{\mathcal{M}}_n }\Biggl|\frac
{1}{n}\sum_{i=1}^n K_h\bigl(r_1(S_i)-x\bigr)\varepsilon_i^* \nonumber
\\[-8pt]
\\[-8pt]
\qquad&& \hspace*{96pt}{} - \frac{1}{n}\sum_{i=1}^n K_h\bigl(r_2(S_i)-x\bigr)\varepsilon_i^* \Biggr| >
Cn^{-\kappa_1} \Biggr)\leq\exp(-cn^c).\nonumber
\end{eqnarray}
Now for $C_I>0$ choose a grid $I_{R,n}$ of $I_R$ with
$O(n^{C_I})$ points, such that for each $x\in I_R$ there exists
a grid point $x^* = x^*(x)\in I_{R,n}$ such that $\|x -
x^*\|\leq n^{-c C_I}$. If $C_I$ is chosen large enough, this
implies that
\begin{equation}\label{t8}
\sup_{x\in I_R}\sup_{r \in\bar{\mathcal{M}}_n }\Biggl|\frac{1}{n}\sum
_{i=1}^n K_h\bigl(r(S_i)-x\bigr)\varepsilon_i
- \frac{1}{n}\sum_{i=1}^n K_h\bigl(r(S_i)-x^*\bigr)\varepsilon_i \Biggr| \leq
n^{-\kappa_1}\hspace*{-45pt}
\end{equation}
for large enough $n$, with probability tending to one.
Furthermore, it follows from \eqref{t7} that
\begin{eqnarray}\label{t9}\quad
\sup_{x\in I_{R,n}}\sup_{r_1,r_2 \in\bar{\mathcal{M}}_n }\Biggl|\frac
{1}{n}\sum_{i=1}^n\! K_h\bigl(r_1(S_i)\!-\!x\bigr)\varepsilon_i
\!-\! \frac{1}{n}\sum_{i=1}^n\! K_h\bigl(r_2(S_i)\!-\!x\bigr)\varepsilon_i \Biggr|\!\leq\!n^{-\kappa_1}. \hspace{-42pt}
\end{eqnarray}
The statement of the lemma then follows
from \eqref{t5}--\eqref{t6} and \eqref{t8}--\eqref{t9}, if
the constants $C_1$ and $C_2$ were chosen large enough.
\end{pf}

\begin{lemma} \label{lemBtheo1}
Suppose that the conditions of Theorem \ref{mainexpansion} hold. Then
\begin{eqnarray*}
&&\sup_{x\in I_R, r_1,r_2\in\bar{\mathcal{M}}_n} \Biggl| \frac{1}{n}\sum
_{i=1}^n K_h\bigl(r_1(S_i)-x\bigr) \biggl(\frac{r_{1,j}(S_i) -x_j}{h_j}\biggr)^a\biggl(\frac
{r_{1,l}(S_i) -x_l}{h_l}\biggr)^b\\
&& \hspace*{60pt}{} -\frac{1}{n}\sum_{i=1}^n K_h\bigl(r_2(S_i)-x\bigr)\biggl(\frac{r_{2,j}(S_i)
-x_j}{h_j}\biggr)^a\biggl(\frac{r_{2,l}(S_i)
-x_l}{h_l}\biggr)^b\Biggr|\\
&&\qquad = O_p\bigl(n^{-(\delta-\eta)_{\mathrm{min}}}\bigr)
\end{eqnarray*}
for $j,l=1,\ldots,q $ $j\neq l$ and $0\leq a+b \leq2$, $0\leq a,b$.
\end{lemma}
\begin{pf}
The lemma follows from
\[
\sup_{x,s}\bigl|K_h\bigl(r_1(s)-x\bigr) - K_h\bigl(r_2(s)-x\bigr) \bigr| \leq C n^{-(\delta-\eta
)_{\mathrm{min}}+\eta_+}
\]
for $r_1,r_2\in\bar{\mathcal{M}}_n$ and from
\begin{eqnarray*}
&&\sup_{x\in I_R,r\in\bar{\mathcal{M}}}\Biggl|\frac{1}{n}\sum_{i=1}^n
K_h\bigl(r(S_i)-x\bigr)\Biggr| \\
&& \qquad  \leq C n^{-1+\eta_+}\sup_{x\in I_R} \#
\{i\dvtx |r_{0,j}(S_i) -x_j| \leq C n^{-\eta_j} \mbox{ for } j=1,\ldots,d\}\\
&& \qquad = O_p(1),
\end{eqnarray*}
which follows from a simple calculation.
\end{pf}

\begin{lemma} \label{lemCtheo1}
Suppose that the assumptions of Theorem \ref{mainexpansion} hold. For a random
variable $R_n = O_p(1)$ that neither depends on $x$ nor $i$, it
holds that
\begin{eqnarray}\qquad
\label{u2}&& \sup_{x\in I_R, 1\leq i \leq n}\bigl|\bigl[m_0(r_0(S_i)) - m_0(x) -
m_0'(x)^T\bigl({r}_0(S_i)-x\bigr)\bigr] I_i(x)\bigr| \nonumber
\\[-8pt]
\\[-8pt]
\qquad&&\qquad \leq R_n n^{-2\eta_{\mathrm{min}}}, \nonumber\\
\label{u3} && \sup_{x\in I_R}\Biggl\|\frac{1}{n}\sum_{i=1}^n K_h\bigl(\hat{r}(S_i)
-x\bigr)\hat w_i(x)\hat w_i(x)^T \nonumber\\
\qquad&& \qquad {} - \frac{1}{n}\sum_{i=1}^n K_h\bigl(r_0(S_i) -x\bigr)\tilde w_i(x)\tilde
w_i(x)^T\Biggr \|  \\
\qquad&&\qquad \leq R_n n^{-(\delta-\eta)_{\mathrm{min}}}, \nonumber\\
\label{u4} && \sup_{x\in I_R}\Biggl\|\frac{1}{n}\sum_{i=1}^n K_h\bigl(r_0(S_i) -x\bigr)\tilde
w_i(x)\tilde w_i(x)^T - f_R(x) B_K \Biggr\|\nonumber
\\[-8pt]
\\[-8pt]
\qquad&&\qquad \leq R_n \bigl(n^{-\eta
_{\mathrm{min}}}+n^{-(1-\eta_+)/2}\sqrt{\log n}\bigr),\nonumber
\end{eqnarray}
where $I_i(x) =\mathbb{I}\{\|(\hat{r}(S_i) -x)/h\|_{1}\leq1\}$ is an
equals one if $\hat{r}(S_i) -x$ lies in
the support of the kernel function $K_h$ and zero otherwise,
and $B_K=\diag(1,\int u^2K(u)\,du,\ldots,\int u^2K(u)\,du)$ is a
$(d+1)\times(d+1)$ diagonal matrix.
\end{lemma}

\begin{pf} Claim (\ref{u2}) follows by a simple calculation. Claim
(\ref{u3}) is a~direct consequence of Lemma \ref{lemBtheo1}, and
(\ref{u4}) follows from standard arguments from kernel
smoothing theory. For the stochastic part, one makes use of
Lemma \ref{genkern}, given in Appendix \ref{genkern_sec}, below.
\end{pf}

\begin{lemma}\label{lemEtheo1} Suppose that the assumptions of Theorem
\ref{mainexpansion} hold. Then
it holds that
\begin{eqnarray}
\label{theo7add1}&\displaystyle\sup_{x\in I_R,r_1,r_2\in\bar{\mathcal{M}}} \| \mu(x, r_1 )- \mu
(x, r_2 )-\E[\mu(x, r_1 )- \mu(x, r_2 )]\|
= O_p(n^{-\kappa_1}),\hspace*{-38pt}&\\
\label{theo7add4}&\displaystyle\sup_{x\in I_R} | \hat\mu(x ) |
= O_p\bigl(\sqrt{\log n}n^{-(1-\eta_+)/2}\bigr),&
\end{eqnarray}
where
\[\hat\mu(x) = n^{-1} \sum_{i=1}^n K_h\bigl(
\hat r(S_i) -x\bigr) \hat w_i(x) \rho(S_i)
\]
 and
 \[\mu(x)=n^{-1}
\sum_{i=1}^n
K_h\bigl( r_0(S_i) -x\bigr)w_i(x) \rho(S_i).
\]
\end{lemma}

\begin{pf} For a proof of (\ref{theo7add1}) one proceeds as in
Lemma \ref{lemAtheo1}.
Claim~(\ref{theo7add4}) follows by classical smoothing arguments. Note
that we have
that $\E( \hat\mu(x, \allowbreak r_0 ))=0$.
\end{pf}

\subsection{\texorpdfstring{Proof of Proposition \protect\ref{proptwostage}}{Proof of Proposition 1}}

In order to prove Proposition \ref{proptwostage}, we use the fact that
the local
polynomial estimator satisfies a certain uniform stochastic expansion
if Assumption
\ref{as4} holds. In order to present this result, we first have to introduce a
substantial
amount of further notation. For simplicity we assume $g_1=\cdots=g_p$, and
we write $g$
for this joint value and for the vector $g=(g,\ldots,g)$.\vadjust{\goodbreak}

Let $N_i = {i+q-1 \choose q-1}$ be the number of distinct $q$-tuples
$u$ with $u_+ =
i$. Arrange these $q$-tuples as a sequence in a lexicographical order
(with the
highest priority given to the last position so that $(0,\ldots, 0, i)$
is the first
element in the sequence, and $(i, 0,\ldots, 0)$ the last element). Let
$\tau_i$
denote this one-to-one mapping, that is, $\tau_i(1) = (0,\ldots, 0,
i)$, \ldots,
$\tau_i(N_i) = (i,0\ldots, 0)$. For each $i=1,\ldots,q$, define a
$N_i\times1$
vector $\mu_i(x)$ with its $k$th element given by $x^{\tau_i(k)}$,
and write
$\mu(x) = (1,\mu_1(x)^T,\ldots,\mu_q(x)^T)^T$, which is a column
vector of length
$N = \sum_{i=1}^q N_i$. Let $\nu_i = \int L(u)u^i\,du$ and define $\nu
_{ni}(x) =
\int L(u)u^if_S(x + gu)\,du$. For $0 \leq j, k \leq q$, let $M_{j,k}$ and
$M_{n,j,k}(x)$ be two $N_j \times N_k$ matrices with their $(l,m)$ elements,
respectively, given by
\[
[ M_{j,k} ]_{l,m} = \nu_{\tau_j(l) + \tau_k(m)} \quad\mbox{and}\quad[
M_{nj,k}(x) ]_{l,m} = \nu_{n,\tau_j(l) + \tau_k(m)}(x).
\]
Now define the $N \times N$ matrices $M_q$ and $M_{n,q}(x)$ by
\begin{eqnarray*}
M_q &=&\pmatrix{M_{0,0} & M_{0,1} &\cdots& M_{0,q}\cr
M_{1,0} & M_{1,1} &\cdots& M_{1,q}\cr
\vdots&\vdots& \ddots&\vdots \cr
M_{q,0} & M_{q,1} &\cdots& M_{q,q}}, \\
M_{n,q}(x) &=&\pmatrix{M_{n,0,0}(x) & M_{n,0,1}(x) &\cdots& M_{n,0,q}(x)\cr
M_{n,1,0}(x) & M_{n,1,1}(x) &\cdots& M_{n,1,q}(x)\cr
\vdots&\vdots& \ddots& \vdots\cr
M_{n,q,0}(x) & M_{n,q,1}(x) &\cdots& M_{n,q,q}(x)}.
\end{eqnarray*}
Finally, denote the first unit $q$-vector by $e_1 = (1,0,\ldots,0)$.
With this
notation, it can be shown along classical lines that the local
polynomial estimator
$\hat{r}$ admits the following stochastic expansion:
\begin{eqnarray}\label{bahadur}
&& \hat{r}(s) = r_0(s) + \frac{1}{n}\sum_{i=1}^n e_1 M_{nq}^{-1}(s)\mu
\bigl((S_i-s)/g\bigr)L_g(S_i-s) \zeta_i \nonumber
\\[-8pt]
\\[-8pt]
&& \hspace*{28pt} {} + g^{q+1}B_n(s) + R_n(s),\nonumber
\end{eqnarray}
where $\sup_{s\in I_S}\| R_n(s)\| = O_p((\log(n)/ng^p)^{1/2})$, and
$B_n$ is a bias term that satisfies
\begin{equation}\quad
B_n(s) = \frac{1}{(q+1)!}e_1 M_{q}^{-1} A_q r_0^{(q+1)}(s) + o_p(1)
\equiv b(s) + o_p(1).\label{fsbias}
\end{equation}
To prove the proposition, define the stochastic component and the bias~term of the
expansion \eqref{bahadur} as \mbox{$\hat r_A(s)\,{=}\, n^{-1}\sum_{i=1}^n\! e_1
M_{nq}^{-1}(s)\mu((S_i\!-\!s)/g)L_g(S_i\!-\!s) \zeta_i $} and $ \hat
r_B(s) \,{=}\, g^{q+1} B_n(s)$, respectively. Now the function $\hat\Delta$ can be
written~as\looseness=-1
\begin{eqnarray*}
\hat\Delta(x) &=& e_1^T N_h(x)^{-1}\E\bigl(K_h\bigl(r_0(S) -x\bigr) w(x,r)\hat
r_A(S)\bigr) \\
&&{} + e_1^T N_h(x)^{-1}\E\bigl(K_h\bigl(r_0(S) -x\bigr) w(x,r)\hat r_B(S)\bigr) +
O_p\biggl(\frac{\log(n)}{ng^p}\biggr)\\
&\equiv&\hat\Delta_A(x) + \hat\Delta_B(x) + O_p\biggl(\frac{\log(n)}{ng^p}\biggr),
\end{eqnarray*}\looseness=0
uniformly over $x\in I_{R}$. We first analyze the term $\hat\Delta_B(x)$.
Through the usual arguments from kernel smoothing theory, one can show
for $x\in
I_{R,n}^-$ that
\begin{eqnarray*}
\hat\Delta_B(x) &=& g^{q+1} e_1^T N_h(x)^{-1}\E\bigl(K_h\bigl(r_0(S) -x\bigr) w(x,r)
b(S)\bigr) + o_p(g^{q+1}) \\
&=& g^{q+1} \E\bigl(b(S)|r_0(S) = x\bigr) + o_p(g^{q+1}+n^{-2\eta})
\end{eqnarray*}
since the function $\E(b(S)|r_0(S) = x)$ is continuous with respect to
$x$ because
of Assumptions \ref{as5} and \ref{as6}. Explicitly, we have
\begin{eqnarray*}
&&\E\bigl(b(S)|r_0(S) = x\bigr) \\
&&\qquad= \frac{ \int b(s_{-p},\varphi(s_{-p},x))
f_S(s_{-p},\varphi(s_{-p},x))\,\partial_{s_{-p}} \varphi(s_{-p},x)
\,ds_{-p} }{ \int f_S(s_{-p},\varphi(s_{-p},x))\,\partial_{s_{-p}}
\varphi(s_{-p},x)\,ds_{-p}}.
\end{eqnarray*}
Next, consider the term $\hat\Delta_A(x)$. Note that for $x\in
I_{R,n}^-$ we have that
\begin{equation}
\hat\Delta_A(x) 
 = \frac{1}{n f_R(x)}\sum_{j=1}^{n} \psi_n(x,S_j) \zeta_j \label{term1}
\end{equation}
with
\begin{eqnarray*}
\psi_n(x,s) &=& \int_{I_S} \bigl( K_h\bigl(r_0(u) -x\bigr) e_1 \bar
M^{-1}_{nq}(u) \mu\bigl((s-u)/g\bigr)L_g(s-u) \bigr) f_S(u)\,du \\
 &=& \int K_h\bigl(r_0(u)
-x\bigr)L_{n,g}^*(s,u - s)\,du,
\end{eqnarray*}
where
$L_{n,g}^*(s,t) = f_S(s-t) e_1 \bar M^{-1}_{nq}(s-t) \mu(t/g)L_g(t)$.
Define $I_{S,n}^-$ as the set that contains all $s \in I_S$ that do not
lie in a
$g$-neighborhood of the boundary of $I_{S}$. Uniformly over $s \in
I_{S,n}^-$, we
have that $M_{n,q}(s)-f_S(s)M_q= O(g)$ . Thus for $ s \in I_{S,n}^-$,
we have that
$\psi_n(x,s) = (1+O(g)) \psi(x,s)$ where the function $\psi$ is
equal to
$\psi(x,s) = \int K_h(r_0(u) -x)L^*_g(u - s)\,du$
with modified kernel $L^*$ defined as
\begin{equation}\label{equivkernel}
L^*(t) = e_1 M_{q}^{-1}\mu(t)L(t).
\end{equation}
Note that $L^*$ is the \emph{equivalent kernel} of the local
polynomial regression
estimator; see \citet{fan1996local}, Section 3.2.2. For $q=0,1$ the equivalent
kernel is in fact equal to the original one, whereas $L^*(t)$ is equal
to $L(t)$
times a polynomial in $t$ of order $q$ for $q\geq2$, with coefficients
such that
its moments up to the order $q$ are equal to zero. The kernel
$L_{n,g}^*(u,t)$ has
the same moment conditions in $t$ as $L^*_g$ but depends on $u$.

We now derive explicit expressions for the leading term in equation
\eqref{term1}
for the cases (a)--(c) of the proposition. Starting with case (a), in
which $g/h \to
0$, it\vadjust{\goodbreak} follows by substitution and Taylor expansion arguments that with
$K_h^{\prime} (v) = h^{-1} K^{\prime} (h^{-1} v)$ and $K_h^{\prime
\prime} (v) =
h^{-1} K^{\prime\prime} (h^{-1} v)$
\begin{eqnarray*}
\psi_n(x,v) &=& \int K_h\bigl(r_0(s) -x\bigr)L^*_{n,g}(s,s - v)\,ds \\
&=& \int K_h\bigl(r_0(v -tg) -x\bigr)L_{n}^*(v-tg,t)\,dt\\
&=& \int\biggl( K_h\bigl(r_0(v)-x\bigr) + K_h'\bigl(r_0(v)-x\bigr) \frac{r_0(v-tg) - r_0(v)}{h}
\\
&&\hspace*{64pt}{}+ K_h''(\chi_1-x)\frac{1}{2}\biggl(\frac{r_0(v-tg) -
r_0(v)}{h}\biggr)^2\biggr)\\
&&\hphantom{\int} {}\times  L_n^*(v-tg,t)\,dt\\
&=& K_h\bigl(r_0(v)-x\bigr)\\
&&{} + K_h'\bigl(r_0(v)-x\bigr) \int\biggl(-\partial_s r_0(v)\frac
{tg}{h} + \partial_s^2 r_0(\chi_2)\frac
{t^2g^2}{2h}\biggr)L_n^*(v-tg,t)\,dt\\
&& {}- \int K_h''(\chi_1-x)\frac{1}{2}\biggl(\frac{\partial_s r_0(\chi
_3)tg}{h}\biggr)^2L_n^*(v-tg,t)\,dt,
\end{eqnarray*}
where $\chi_1$, $\chi_2$ and $\chi_3$ are intermediate values between
$r_0(v)$ and
$r_0(v -tg)$, $v$ and $v-tg$, and $v$ and $v-tg$, respectively. This
gives an
expansion for $\psi_n(x,v)$ of order $(g/h)^2$. For $v \notin
I_{S,n}^-$ one gets
an expansion of order $g/h$. Put $k_n(v) = - \partial_s r_0(v) \int t
L_n^*(v-tg,t)\,dt$.
Together with Lemma \ref{genkern} in Appendix~\ref{genkern_sec}, we
thus obtain
that\vspace*{-1pt}
\begin{eqnarray*}
&&\frac{1}{n f_R(x)}\sum_{j=1}^{n} \psi_n(x,S_j) \zeta_j \\
&&\qquad =\frac{1}{n f_R(x)}\sum_{i=1}^n \biggl(K_h\bigl(r_0(S_i)-x\bigr)+ \frac{g}{h}K^{\prime}_h
\bigl(r_0(S_i)-x\bigr) k_n(S_i)\biggr) \zeta_i\\
&& \hspace*{30pt} {}+ O_p\biggl(\biggl(\frac{g}{h}\biggr)^2\biggl(\frac{\log(n)}{nh}\biggr)^{1/2}\biggr)\\
&&\qquad =\frac{1}{n f_R(x)}\sum_{i=1}^n K_h\bigl(r_0(S_i)-x\bigr)\zeta_i + O_p\Biggl(\Biggl(\frac
{g^2}{h^2}+ \sqrt{\frac{g^3}{h^2}}\Biggr)\sqrt{\frac{\log(n)}{nh}}\Biggr),
\end{eqnarray*}
as claimed. To show statement (b) of the proposition, we rewrite the function~$\psi_n$ as follows:\vspace*{-1pt}
\begin{eqnarray*}
\psi_n(x,v) &=& \int\biggl(K_h\bigl(r_0(v) -x + \partial_s r_0(v)th \bigr) + K'\biggl(\frac
{\chi_1}{h}\biggr)\partial_s^2r_0(\chi_2)\frac{1}{2}t^2\biggr)\\
&& \hspace*{7pt} {}\times L_n^*(v-th,t)\,dt\\
&=& J_{n,h}(x,v) + h \int K_h'(\chi_1)\partial_s^2r_0(\chi_2)\frac
{1}{2}t^2 L_n^*(v-th,t)\,dt,
\end{eqnarray*}
where $J_{n,h}(x,s) = \int K_h(r_0(s)-x - \partial_s r_0 (s) uh)
L_n^*(s-uh,u)\,du$, and $\chi_1$ is an intermediate value between
$r_0(v +gt)$ and
$r_0(v) + \partial_s r_0(v)tg $, and $\chi_2$ is an intermediate
value between $v$
and $v+gt$. As in the proof of part (a), it follows from Lemma \ref{genkern} in
Appendix \ref{genkern_sec} that
\begin{eqnarray*}
\frac{1}{n f_R(x)}\sum_{j=1}^{n} \psi_n(x,S_j) \zeta_j &=& \frac
{1}{n f_R(x)}\sum_{j=1}^{n} J_{n,h}(x,S_j)\zeta_j + O_p\Biggl(h \sqrt
{\frac{\log(n)}{nh}}\Biggr)\\
&=& \frac{1}{n f_R(x)}\sum_{j=1}^{n} J_{h}(x,S_j)\zeta_j + O_p\Biggl(\sqrt
{\frac{\log(n)}{n}}\Biggr),
\end{eqnarray*}
where $J_{h}$ uses the location independent form of the equivalent
kernel $L^*$ as
defined in the text in front of Proposition \ref{proptwostage}. This implies the desired result.

Now consider statement (c) of
the proposition. In this case, where $g/h \to\infty$, we can rewrite
the function
$\psi_n$ as follows:
\begin{eqnarray*}
&&\psi_n(x,v) = \int K_h(w_p -x)\\
&& \hspace*{58pt} {} \times L_{n,g}^*\bigl(\bigl(w_{-p}, \varphi(w)\bigr )^T,
\bigl(w_{-p} - v_{-p}, \varphi(w) - v_p\bigr)^T\bigr)\,\partial_x \varphi(w)\,dw.
\end{eqnarray*}
From tedious but conceptually simple Taylor expansion arguments similar
to the
ones employed for case (a), and from Lemma \ref{genkern}, one gets that
\[
\frac{1}{n f_R(x)}\sum_{j=1}^{n} \psi_n(x,S_j) \zeta_j = \frac
{1}{n f_R(x)}\sum_{j=1}^{n} H_{n,g}(x,S_j)\zeta_j + O_p\Biggl(\frac
{h^2}{g^2} \sqrt{\frac{\log(n)}{ng}}\Biggr),
\]
where
\begin{eqnarray}\label{defHn}
&&H_{n,g}(x,v) = \int K(t) L_{n,g}^*\bigl( \bigl(v_{-p} + g s_{-p}, G_n(v_{-p},
x;s_{-p}, t)\bigr),\nonumber \\
&& \hspace*{132pt}\bigl(s_{-p},G_n(v_{-p}, x;s_{-p}, t)-v_p\bigr) \bigr)  \\
&& \hspace*{70pt}\partial_x \varphi(v_{-p},x)\,ds_{-p}\,dt\nonumber
\end{eqnarray}
and $G_n(v_{-p}, x; s_{-p}, t) = \varphi(v_{-p},x) + g
s_{-p}\,\partial_{{-p}}\varphi(v_{-p},x) + ht\,\partial_x \varphi
(v_{-p},x)$. With $H^{\Delta}_n$ as defined in the text, we find
\begin{eqnarray*}
\frac{1}{n f_R(x)}\sum_{j=1}^{n} \psi_n(x,S_j) \zeta_j &=& \frac
{1}{n f_R(x)}\sum_{j=1}^{n} H^{\Delta}_n(x,S_j)\zeta_j \\
&& {} + O_p\Biggl(\Biggl(1+\sqrt{\frac{h}{g}}\Biggr) \sqrt{\frac{\log(n)}{n}}+\frac
{h^2}{g^2} \sqrt{\frac{\log(n)}{ng}}\Biggr).
\end{eqnarray*}
Since $O(h/g)=o(1)$, this completes our proof.

\subsection{\texorpdfstring{Proof of Proposition \protect\ref{proptwostage2b}}{Proof of Proposition 2}}
To show the result, note that
\begin{eqnarray*}
\Gamma(x,r) &=& e_1^T N_h(x)^{-1}\E
\bigl(\bigl(K_h\bigl(r(S)-x\bigr)-K_h\bigl(r_0(S)-x\bigr)\bigr)w(x)\rho(S)\bigr) \\
&&{}+ O_p\bigl(n^{-((1/2)(1-\eta_{+})+2\delta-\eta)}\bigr)\\
&=& \E\bigl(\rho(S)|r(S)=x\bigr) - \E\bigl(\rho(S)|r_0(S)=x\bigr) \\
&&{}+ O_p\bigl(n^{-2\eta} +
n^{-((1/2)(1-\eta_{+})+2\delta-\eta)}\bigr)
\end{eqnarray*}
uniformly over $x\in I_R$ and $r\in\mathcal{M}_n$. Since $\E(\rho
(S)|r_0(S))\equiv0$
by construction, it suffices to consider the term $\E(\rho
(S)|r(S)=x)$. To simplify the
exposition, we strengthen Assumption \ref{as6} and suppose that in addition to
$r_0$ all functions
$r\in\mathcal{M}_n$ are strictly monotone with respect to their last
argument, and write~$\varphi_r$ for corresponding the inverse function that satisfies
$r(u_{-p},\varphi_r(u_{-p},x))=x$
(without this condition, the notation would be much more involved, as
we would have
to consider all regions where the functions $r\in\mathcal{M}_n$ are piecewise
monotone with respect to the last component separately). Using rules
for integrals on manifolds, we
derive the following explicit expression for $\E(\rho(S)|r(S)=x)$:
\begin{eqnarray*}
&&\E\bigl(\rho(S)| r(S) =x\bigr) \\
&& \qquad = \frac{\int\rho(s_{-p},\varphi
_r(s_{-p},x))f_S(s_{-p},\varphi_r(s_{-p},x))\,\partial_{-p} \varphi
_r(s_{-p},x)\,ds_{-p}}{\int f_S(s_{-p},\varphi_r(s_{-p},x))\,\partial
_{-p} \varphi_r(s_{-p},x)\,ds_{-p}}.
\end{eqnarray*}
Set the numerator of the above expression as $\gamma_1(x,r)$ and the
denominator as $\gamma_2(x,r)$.
Then clearly $\gamma_2(x,\hat r)=f_R(x) + o_p(1)$ uniformly over
$x\in I_R$. Moreover,
note that the mapping
\[
r \mapsto\rho(s_{-p},\varphi_r(s_{-p},x))f_S(s_{-p},\varphi_r(s_{-p},x))
\]
is Hadamard differentiable at $r_0$, with derivative
\[
r\mapsto\frac{\partial_p \lambda(s_{-p},\varphi
(s_{-p},x))}{\partial_{p}r_{0}(s_{-p},\varphi(s_{-p},x)) }
r(s_{-p},\varphi(s_{-p},x)).
\]
It follows with $\gamma_1(x,r_0)=0$ that
\begin{eqnarray*}
\gamma_1(x,r)&=& \int\frac{\partial_p \lambda(s_{-p},\varphi
(s_{-p},x))}{\partial_{p}r_{0}(s_{-p},\varphi(s_{-p},x)) }
\bigl(r(s_{-p},\varphi(s_{-p},x))-r_0(s_{-p},\varphi(s_{-p},x))\bigr) \\
&&\hspace{7pt}{}\times(\partial_{-p} \varphi_r(s_{-p},x))\,ds_{-p}\\
 &&\hspace{-4pt} {}+ O_p(\| r -r_0\|_\infty^2).
\end{eqnarray*}
We evaluate the term $\gamma_1(x,\hat r)$, substitute the uniform
expansion \eqref{bahadur} for $\hat r(s) - r_0(s)$ into the explicit expression
derived above, and use standard arguments from kernel smoothing theory.
This gives
the desired expansion for $\hat\Gamma_A$. The form of $\hat\Gamma
_B$ follows from
the same arguments used to derive the form of $\hat\Delta_B$ in the
proof of
Proposition \ref{proptwostage}.

\subsection{\texorpdfstring{Proofs of Corollaries \protect\ref{corAtwostage}--\protect\ref{corCtwostage}}
{Proofs of Corollaries 1--4}}
The statements of these corollaries follow by direct
application of Proposition \ref{proptwostage}--\ref{proptwostage2b} and Theorem \ref{mainexpansion}. The
statement of Corollary \ref{corAtwostage} is immediate. For
Corollaries \ref{corBtwostage}--\ref{corCtwostage},
we only have to check that the error bounds in Theorem \ref{mainexpansion} and
Proposition \ref{proptwostage}--\ref{proptwostage2b} are of the desired order. We
only discuss how the constants $\alpha$, $\delta$ and $\xi$ can
be chosen. Note that all these constants have no subindex
because we only consider the case $d=1$. We apply Theorem \ref{mainexpansion}
conditionally on the values of $S_1,\ldots,S_n$. Then the only
randomness in the pilot estimation comes from
$\zeta_1,\ldots,\zeta_n$. We can decompose $\hat r$ into
$\hat r_A+ \hat r_B$, where $\hat r_A$ is the local
polynomial fit to $(S_i,\zeta_i)$, and~$\hat r_B$ is the
local polynomial fit to $(S_i,r_0(S_i))$. Conditionally given
$S_1,\ldots,S_n$, the value of $\hat r_B$ is fixed, and for
checking Assumption~\ref{as3}, we only have to consider entropy
conditions for sets of possible outcomes of $\hat r_A$. We
will show that with $\alpha= p/k$ one can choose for $\delta$
and $\xi$ any value that is larger than $(1-p \theta) /2$ or
$-pk^{-1}(1-p \theta) /2 +p \theta$, respectively. Note that
then $\alpha\leq2$ because of Assumption \ref{as4}(iii). It can be
easily checked that we get the desired expansions in
Corollaries \ref{corAtwostage} and \ref{corBtwostage} with this
choices of $\alpha= p/k$,~$\delta$ and~$\xi$ (with~$\delta$
and~$\xi$ small enough). In particular note that we can make
$\delta\alpha+\xi$ as close to $p \theta$ as we like.

It is clear that Assumption \ref{as2} holds for this choice of
$\delta$. This follows by standard smoothing theory for local
polynomials. Compare also Lemma \ref{genkern}
and the proof of Proposition \ref{proptwostage}. It remains to check
Assumption \ref{as3}.
It suffices to check the entropy conditions for the tuple of
functions $(n^{-1} \sum_{i=1}^n L_h(S_i-s)
[(S_i-s)/g]^{\pi}\zeta_i\dvtx 0\leq\pi_+ \leq q, \pi_j\geq0$ for
$j=1,\ldots,p)$. This follows because we get $\hat r_A$ by
multiplying this tuple of functions with a (stochastically)
bounded vector. We now argue that all derivatives of order $k$
of the functions $n^{-1} \sum_{i=1}^n L_h(S_i-s)
[(S_i-s)/g]^{\pi}\zeta_i$ can be bounded by a variable $B_n$
that fulfills $B_n \leq b_n =n^{\xi^{**}})$ with probability
tending to one. Here $\xi^{**}$ is a number with $\xi^{**} > -
\frac{1 }{ 2} (1- p\theta)+k\theta$. This bound holds uniformly
in $s$ and $\pi$. Furthermore, the functions $n^{-1}
\sum_{i=1}^n L_h(S_i-s) [(S_i-s)/g]^{\pi}\zeta_i$ can be
bounded by a variable $A_n$ that fulfills $A_n \leq a_n
=n^{\xi^*})$ with probability tending to one. Here $\xi^*$ is a
number with $\xi^* > - \frac{1 }{ 2} (1- p\theta)$. Again, this
bound holds uniformly in $s$ and $\pi$. We now consider the set
of functions on $I_S$ that are absolutely bounded by $a_n$ and
that have all partial derivatives of order $k$ absolutely
bounded by $b_n$. We argue that this set can be covered by $C
\exp( \lambda^{-p/k} b_n ^{p/k}) $ balls with $\|\cdot
\|_{\infty}$-radius $\lambda$ for $\lambda\leq a_n$. Here the
constant $C$ does not depend on $a_n$ and $b_n$. This entropy
bound shows that Assumption \ref{as3} holds with these choices of
$\alpha$, $\delta$ and $\xi$. For the proof of the entropy
bound one applies an entropy bound for the set of functions on
$I_S$ that are absolutely bounded by $1$ and that have all
partial derivatives of order $k$ absolutely bounded by $1$.
This set can be covered by $C \exp( \lambda^{-p/k}) $ balls
with $\|\cdot\|_{\infty}$-radius $\lambda$ for $\lambda\leq
1$. The desired entropy bound follows by rescaling of the
functions. Note that we have that $b_n^{-1}a_n \to0$.

\subsection{\texorpdfstring{Proof of Corollary \protect\ref{corcens}}{Proof of Corollary 5}}
Our proof has the same structure as the one provided by
\citet{linton2002nonparametric}, but making use of Theorem \ref{mainexpansion}
considerably simplifies some of their arguments. First, note
that the restriction that $\underline{\theta}<\theta<
\bar{\theta}$ implies that $(ng^p)^{1/2}h^2\to0$ and
$(ng^p)^{1/2}g^{q+1}\to0$. From a~second-order Taylor
expansion, we furthermore obtain that
\begin{eqnarray*}
\hat{\mu}(x)-\mu_0(x)&=&\frac{1}{q_0(r_0(x))}\bigl(\hat{r}(x)-r_0(x)\bigr)\\[-2pt]
&&{} + \int_{r_0(x)}^{\lambda}\frac{\hat{q}(s)-q_0(s)}{q_0(s)^2}\,ds-
 \frac{\hat{q}'(\bar{r}(x))}{2\hat{q}(\bar{r}(x))^2}\bigl(\hat
{r}(x)-r(x)\bigr)^2 \\[-2pt]
&&{} - \int_{r(x)}^{\lambda}\frac{(\hat{q}(s)-q_0(s))^2}{\hat
{q}(s)q_0(s)^2}\,ds\\[-2pt]
&&{} + \frac{(\hat{q}(\check{r}(x))-q_0(\check{r}(x)))^2}{\hat{q}(\check{r}(x))q_0(\check
{r}(x))}\bigl(\hat{r}(x)-r_0(x)\bigr)\\[-2pt]
&\equiv& T_1 + T_2 + T_3 + T_4 + T_5,
\end{eqnarray*}
where $\hat{r}(x)$ and $\check{r}(x)$ are intermediate values
between $r(x)$ and $\hat{r}(x)$. Now it follows from standard
arguments for local linear estimators that
\[
\sqrt{ng^p}T_1 \stackrel{d}{\rightarrow} N\biggl(0,\frac{\sigma
^2_r(x)}{f_S(x)s_0^2(x)}\int L^2(t)\,dt\biggr),
\]
since $s_0(x) = q_0(r_0(x))$. To prove the corollary, it thus
only remains to be shown that the remaining four terms in the
above expansion are of smaller order than $T_1$. Under the
conditions of the corollary, it is easy to show with
straightforward rough arguments that
$\inf q(s) > 0$, $\sup\hat{q}'(s) = O_p(1)$ and $\sup|\hat
{q}(s)-q_0(s)|^2 = o_p((ng^p)^{-1/2})$
where the supremum and infimum
are taken over $ s\in(r_o(x)-\epsilon,\lambda_0 + \epsilon)$
for some $\epsilon> 0$, respectively. 
This directly implies that $T_3+T_4 + T_5 =
o_p((ng^p)^{-1/2})$. Now consider the term $T_2$. From Theorem
\ref{mainexpansion}, we obtain that
\[
T_2 = \int_{r_0(x)}^{\lambda}\frac{\tilde
{q}(s)-q_0(s)}{q_0(s)^2}\,ds -\int_{r_0(x)}^{\lambda}\frac{q_0'(s)\hat
\Delta(s) - \hat\Gamma(s)}{q_0(s)^2}\,ds + O_p(n^{-\kappa}),
\]
where $\tilde{q}(x)$ is the oracle estimator of the function $q$
obtained via local
linear regression of $\mathbb{I}\{Y>0\}$ on $r_0(X)$, and $\hat\Delta
(s)$ and
$\hat\Gamma(x)$
are the adjustment terms that appear in the main expansion in Theorem
\ref{mainexpansion}, with the necessary
adjustments to the notation. Using similar
arguments as in the proof of Proposition \ref{proptwostage}--\ref{proptwostage2b} and
Corollaries \ref{corBtwostage}--\ref{corCtwostage},
and the restriction that $\underline{\theta}<\theta<\bar{\theta}$,
we obtain that
\begin{eqnarray*}
\int_{r(x)}^{\lambda}\frac{\tilde{q}(s)-q(s)}{q^2(s)}\,ds
&=& \frac
{1}{n}\sum_{i=1}^n \frac{\varepsilon_i}{f_R(r_0(X_i))} + O_p(h^2) \\[-2pt]
&=& O_p(n^{-1/2})+O_p(h^2) = o_p((ng^p)^{-1/2})
\end{eqnarray*}
for $\varepsilon_i = \mathbb{I}\{Y_i > 0\} - q_0(X_i)$, and similarly that
\begin{eqnarray*}
\int_{r(x)}^{\lambda}\frac{q_0'(s)\hat\Delta(s) - \hat\Gamma
(s)}{q_0(s)^2}\,ds &=& O_p(n^{-1/2}) +
O_p\biggl(\frac{\log{n}}{ng^p} \biggr) + O_p(g^{q+1}) \\
&=& o_p((ng^p)^{-1/2}).
\end{eqnarray*}
Thus $T_2= o_p((ng^p)^{-1/2})$.
Finally, straightforward calculations show that
$\underline{\theta}<\theta<\bar{\theta}$ also implies that
$O_p(n^{-\kappa}) = o_p((ng^p)^{-1/2})$. This completes the
proof.

\subsection{\texorpdfstring{Proof of Corollary \protect\ref{corsim}}{Proof of Corollary 6}}

Let $\hat{f} = (\hat{m},\hat{\mu}_2)$ and $\bar{f} =
(m,\mu_2)$, define the functional $S_n(f)$ as
\[
S_n(f) = \frac{1}{n}\sum_{i=1}^n f_1\bigl(x_1,z_1,X_{1i} - f_2(Z_i)\bigr) - \mu
_1(x_1,z_1),
\]
and let $\dot{S}_n(f)[h] =\lim_{t \to0} (S_n(f+th)-S_n(f))/t$
denote its directional derivative. One then obtains through
direct calculations that for any $f = (f_{1,A}+f_{1,B},f_2)$
with bounded second derivatives we have that
\begin{eqnarray*}
&&\|S_n(f)-S_n(\bar{f}) - \dot{S}_n(\bar{f})[f-\bar{f}] \|_\infty\\
&&\qquad= O( \|f_2-\bar{f}_{2}\|_\infty^2) + O\bigl(\|f_2-\bar{f}_{2}\|_\infty\bigl\|
f_{1,A}^{(v)}-\bar{f}_{1}^{(v)}\bigr\|_\infty\bigr) + O(\|f_{1,B}\|_\infty),
\end{eqnarray*}
where $f_{1,A}^{(v)}(x_1,z_1,v) = \partial_v f_{1,A}(x_1,z_1,v)$.
Using the same kind of arguments as in the proof of Proposition
\ref{proptwostage}, under the conditions of the corollary one
can derive the following stochastic expansion of
$\hat{m}$ up to order $o_p((nh^{1+d_1})^{-1/2})$, uniformly over $(x_1,z_1,v)$
in the $h$-interior of the support of $(X_1,Z_1,V)$:
\begin{eqnarray}\label{expsi}
&&\hat{m}(x_1,z_1,v) - m(x_1,z_1,v) \nonumber\\
&& \qquad = \frac{1}{n f_R(x_1,z_1,v)}\sum_{i=1}^n K_h
\bigl((X_{1i},Z_{1i},V_{i})-(x_1,z_1,v)\bigr)\varepsilon_i \\
&& \qquad \quad {} + o_p((nh^{1+d_1})^{-1/2}), \nonumber
\end{eqnarray}
where $\varepsilon_i = Y - m(X_{1i},Z_{1i},V_{i})$. A similar, but notationally
more involved expansion can be derived for values of $(x_1,z_1,v)$ in
the proximity
of the boundary. Note that since exclusion
restriction on the instruments that $\E(U|Z_1,Z_2,V)=\E(U|V)$ implies that
$\E(\varepsilon|Z_1,Z_2,V)=0$. In the notation of Theorem \ref{mainexpansion}, this
means that
$\rho(s)\equiv0$, and hence the term corresponding to $\hat\Gamma
(x)$ is equal
to zero and does not need to be considered.

Now let $\hat{f}_{1,A}$ denote the sum of the function $m$ and the
leading term of the expansion \eqref{expsi},
and denote the remainder term by $\hat{f}_{1,B}$. Then it follows
from, for example, \citet{masry1996multivariate} and the conditions on
$\eta$ and $\theta$, that
\[
\|\hat{f}_2 -\bar{f}_2\|_\infty=
O_P\bigl(\bigl(\log(n)/(ng^{d_1+d_2})\bigr)^{1/2} \bigr) =
o_p((nh^{1+d_1})^{-{1/4}}),
\]
and it follows from the same
result together with Lemma \ref{genkern}
in Appendix~\ref{genkern_sec} that
\begin{eqnarray*}
\|\hat{f}_2-\bar{f}_{2}\|_\infty\bigl\|\hat{f}_{1,A}^{(v)}-\bar
{f}_{1}^{(v)}\bigr\|_\infty
&=& O_P\bigl(\log(n)/(n^2h^{3+d_1} g^{d_1+d_2})^{1/2} \bigr)\\[-2pt]
&=& o_p((nh^{1+d_1})^{-{1/2}}).
\end{eqnarray*}
For any fixed values $(x_1,z_1)$
we thus have that
\[
\hat{\mu_1}(x_1,z_1)- \mu_1(x_1,z_1) = S_n(\hat{f}) = S_n(\bar{f})
+ T_{1,n} + T_{2,n} + o_p((nh^{1+d_1})^{-{1/2}}),
\]
where
\begin{eqnarray*}
T_{1,n} &=& -\frac{1}{n}\sum_{i=1}^n m^{(v)}(x_1,z_1,V_i)\bigl(\hat{\mu
}_2(Z_i) - \mu_2(Z_i)\bigr),\\[-2pt]
T_{2,n} &=& \frac{1}{n}\sum_{i=1}^n \bigl(\hat{m}(x_1,z_1,V_i) - m(x_1,z_1,V_i)\bigr).
\end{eqnarray*}
Being a simple sample average of i.i.d. mean zero random
variables, one can directly see that $S_n(f_0) =
O_p(n^{-1/2})=o_p((nh^{1+d_1})^{-{1/2}})$. Using a stochastic
expansion for $\hat\mu_2$ as in the proof of Proposition
\ref{proptwostage}, and applying projection arguments for
U-statistics, one also finds that $T_{1,n} =
O_p(n^{-1/2})=o_p((nh^{1+d_1})^{-{1/2}})$. Now consider the
term $T_{2,n}$. From the expansion in \eqref{expsi}, it follows
that for any fixed values $(x_1,z_1)$ we have that
\begin{eqnarray}\label{mu1exp}\qquad
T_{2,n} &=& \frac{1}{n}\sum_{j=1}^n\frac{1}{n f_R(x_1,z_1,V_j)}\sum
_{i=1}^n K_h
\bigl((X_{1i},Z_{1i},V_{i})-(x_1,z_1,V_j)\bigr)\varepsilon_i \nonumber
\\[-9pt]
\\[-9pt]
\qquad&& {}+ o_p((nh^{1+d_1})^{-1/2})\nonumber.
\end{eqnarray}
This in turn implies that
\[
\sqrt{nh^{1+d_1}}T_{2,n} \stackrel{d}{\rightarrow}N\biggl(0, \E\biggl(\frac
{\sigma_\varepsilon^2(x_1,z_1,V)}{f_{XZ_1|V}(x_1,z_1,V)}\biggr) \int\tilde
{K}(t)^2\,dt\biggr)
\]
using again projection arguments for U-statistics.\vspace*{-3pt}

\subsection{Uniform rates for generalized
kernels}\label{genkern_sec}

The following auxiliary lemma states uniform rates for averages of
i.i.d. mean zero random variables weighted by ``kernel-type''
expressions. It is used in the proofs of several of our
results. Modifications of the lemma are well known in the
smoothing literature; see, for example, \citet{haerdle1988uniform}. The
lemma can be proved by standard smoothing arguments. One can
proceed by using a Markov inequality as in the proof of Lemma
\ref{lemAtheo1}, but without making use of a chaining argument.

\begin{lemma}\label{genkern}
Assume that $D\subset\mathbb{R}^{d_x}$ is a compact set, and
$W_{n,h}$ is a kernel-type function that satisfies
$W_{n,h}(u,z)=0$ for $\Vert u-t(z)\Vert> b_n h$ for some deterministic\vadjust{\goodbreak}
sequence $0<b\leq|b_n|\leq B<\infty$, and
$t\dvtx\mathbb{R}^{d_{S}}\to\mathbb{R}^{d_x}$ a continuously
differentiable function, for any $u \in D$ and $z\in
\mathbb{R}^{d_{S}}$. Furthermore, assume that
$|W_{n,h}(u,z)-W_{n,h}(v,z)|\leq l
\frac{\Vert u-t(z)\Vert}{h}h^{-d_x}\tilde{W}_n(v,t(z))$ with
$\sup_n\tilde{W}_n$ bounded, and that
$\mathbb{E}[\exp{(\rho|\varepsilon|)}|S]<C$ a.s. for a constant
$C>0$ and $\rho> 0$ small enough. Then we have that
\[
\sup_{x\in D}\Biggl|\frac{1}{n}\sum_{i=1}^n a_n W_{n,h}(x,S_i)\varepsilon
_i\Biggr| = O_p\Biggl(\sqrt{\frac{\log(n)}{nh^{d_x}}}\Biggr)
\]
for any deterministic sequence $a_n$ with $|a_n|\leq A$.
\end{lemma}
\end{appendix}

\section*{Acknowledgments}
We would like to thank the Associate Editor and three
anonymous referees for
their comments.



%

\printaddresses


\begin{thebibliography}{34}

\bibitem[\protect\citeauthoryear{Ahn}{1995}]{ahn1995nonparametric}
\begin{barticle}[mr]
\bauthor{\bsnm{Ahn},~\bfnm{Hyungtaik}\binits{H.}}
(\byear{1995}).
\btitle{Nonparametric two-stage estimation of conditional choice probabilities
  in a binary choice model under uncertainty}.
\bjournal{J. Econometrics}
\bvolume{67}
\bpages{337--378}.
\bid{doi={10.1016/0304-4076(94)01636-E}, issn={0304-4076}, mr={1333107}}
\bptok{imsref}%
\end{barticle}
\endbibitem

\bibitem[\protect\citeauthoryear{Andrews}{1994}]{andrews1994asymptotics}
\begin{barticle}[mr]
\bauthor{\bsnm{Andrews},~\bfnm{Donald W.~K.}\binits{D.~W.~K.}}
(\byear{1994}).
\btitle{Asymptotics for semiparametric econometric models via stochastic
  equicontinuity}.
\bjournal{Econometrica}
\bvolume{62}
\bpages{43--72}.
\bid{doi={10.2307/2951475}, issn={0012-9682}, mr={1258665}}
\bptok{imsref}%
\end{barticle}
\endbibitem

\bibitem[\protect\citeauthoryear{Andrews}{1995}]{andrews1995nonparametric}
\begin{barticle}[mr]
\bauthor{\bsnm{Andrews},~\bfnm{Donald W.~K.}\binits{D.~W.~K.}}
(\byear{1995}).
\btitle{Nonparametric kernel estimation for semiparametric models}.
\bjournal{Econometric Theory}
\bvolume{11}
\bpages{560--596}.
\bid{doi={10.1017/S0266466600009427}, issn={0266-4666}, mr={1349935}}
\bptok{imsref}%
\end{barticle}
\endbibitem

\bibitem[\protect\citeauthoryear{Blundell and
  Powell}{2004}]{blundell2004endogeneity}
\begin{barticle}[mr]
\bauthor{\bsnm{Blundell},~\bfnm{Richard~W.}\binits{R.~W.}} \AND
  \bauthor{\bsnm{Powell},~\bfnm{James~L.}\binits{J.~L.}}
(\byear{2004}).
\btitle{Endogeneity in semiparametric binary response models}.
\bjournal{Rev. Econom. Stud.}
\bvolume{71}
\bpages{655--679}.
\bid{doi={10.1111/j.1467-937X.2004.00299.x}, issn={0034-6527}, mr={2062893}}
\bptok{imsref}%
\end{barticle}
\endbibitem

\bibitem[\protect\citeauthoryear{Chen, Linton and
  Van~Keilegom}{2003}]{chen2003estimation}
\begin{barticle}[mr]
\bauthor{\bsnm{Chen},~\bfnm{Xiaohong}\binits{X.}},
  \bauthor{\bsnm{Linton},~\bfnm{Oliver}\binits{O.}} \AND
  \bauthor{\bsnm{Van~Keilegom},~\bfnm{Ingrid}\binits{I.}}
(\byear{2003}).
\btitle{Estimation of semiparametric models when the criterion function is not
  smooth}.
\bjournal{Econometrica}
\bvolume{71}
\bpages{1591--1608}.
\bid{doi={10.1111/1468-0262.00461}, issn={0012-9682}, mr={2000259}}
\bptok{imsref}%
\end{barticle}
\endbibitem

\bibitem[\protect\citeauthoryear{Conrad and Mammen}{2009}]{conrad2009garch}
\begin{bmisc}[author]
\bauthor{\bsnm{Conrad},~\bfnm{C.}\binits{C.}} \AND
  \bauthor{\bsnm{Mammen},~\bfnm{E.}\binits{E.}}
(\byear{2009}).
\btitle{{Nonparametric regression on a generated covariate with an application
  to semiparametric GARCH-in-Mean models}}.
\bhowpublished{Unpublished manuscript}.
\bptok{imsref}%
\end{bmisc}
\endbibitem

\bibitem[\protect\citeauthoryear{Das, Newey and
  Vella}{2003}]{das2003nonparametric}
\begin{barticle}[mr]
\bauthor{\bsnm{Das},~\bfnm{Mitali}\binits{M.}},
  \bauthor{\bsnm{Newey},~\bfnm{Whitney~K.}\binits{W.~K.}} \AND
  \bauthor{\bsnm{Vella},~\bfnm{Francis}\binits{F.}}
(\byear{2003}).
\btitle{Nonparametric estimation of sample selection models}.
\bjournal{Rev. Econom. Stud.}
\bvolume{70}
\bpages{33--58}.
\bid{doi={10.1111/1467-937X.00236}, issn={0034-6527}, mr={1952565}}
\bptok{imsref}%
\end{barticle}
\endbibitem

\bibitem[\protect\citeauthoryear{d'Haultfoeuille and
  Maurel}{2009}]{dHautfoeuille2009Roy}
\begin{bmisc}[author]
\bauthor{\bsnm{d'Haultfoeuille},~\bfnm{X.}\binits{X.}} \AND
  \bauthor{\bsnm{Maurel},~\bfnm{A.}\binits{A.}}
(\byear{2009}).
\btitle{{Inference on a generalized Roy model, with an application to schooling
  decisions in France}}.
\bhowpublished{Unpublished manuscript}.
\bptok{imsref}%
\end{bmisc}
\endbibitem

\bibitem[\protect\citeauthoryear{Einmahl and Mason}{2000}]{einmahl2000uniform}
\begin{barticle}[mr]
\bauthor{\bsnm{Einmahl},~\bfnm{Uwe}\binits{U.}} \AND
  \bauthor{\bsnm{Mason},~\bfnm{David~M.}\binits{D.~M.}}
(\byear{2000}).
\btitle{An empirical process approach to the uniform consistency of kernel-type
  function estimators}.
\bjournal{J. Theoret. Probab.}
\bvolume{13}
\bpages{1--37}.
\bid{doi={10.1023/A:1007769924157}, issn={0894-9840}, mr={1744994}}
\bptok{imsref}%
\end{barticle}
\endbibitem

\bibitem[\protect\citeauthoryear{Escanciano, Jacho-Ch{\'a}vez and
  Lewbel}{2011}]{escanciano2011uniform}
\begin{bmisc}[author]
\bauthor{\bsnm{Escanciano},~\bfnm{J.~C.}\binits{J.~C.}},
  \bauthor{\bsnm{Jacho-Ch{\'a}vez},~\bfnm{D.}\binits{D.}} \AND
  \bauthor{\bsnm{Lewbel},~\bfnm{A.}\binits{A.}}
(\byear{2011}).
\btitle{{Uniform convergence for semiparametric two step estimators and
  tests}}.
\bhowpublished{Unpublished manuscript}.
\bptok{imsref}%
\end{bmisc}
\endbibitem

\bibitem[\protect\citeauthoryear{Fan and Gijbels}{1996}]{fan1996local}
\begin{bbook}[author]
\bauthor{\bsnm{Fan},~\bfnm{J.}\binits{J.}} \AND
  \bauthor{\bsnm{Gijbels},~\bfnm{I.}\binits{I.}}
(\byear{1996}).
\btitle{{Local Polynomial Modelling and Its Applications}}.
\bpublisher{CRC Press}, \baddress{New York}.
\bptok{imsref}%
\end{bbook}
\endbibitem

\bibitem[\protect\citeauthoryear{Hahn and Ridder}{2011}]{hahn2010}
\begin{bmisc}[author]
\bauthor{\bsnm{Hahn},~\bfnm{J.}\binits{J.}} \AND
  \bauthor{\bsnm{Ridder},~\bfnm{G.}\binits{G.}}
(\byear{2011}).
\btitle{{The asymptotic variance of semiparametric estimators with generated
  regressors}}.
\bhowpublished{Unpublished manuscript}.
\bptok{imsref}%
\end{bmisc}
\endbibitem

\bibitem[\protect\citeauthoryear{H{\"a}rdle, Janssen and
  Serfling}{1988}]{haerdle1988uniform}
\begin{barticle}[mr]
\bauthor{\bsnm{H{\"a}rdle},~\bfnm{W.}\binits{W.}},
  \bauthor{\bsnm{Janssen},~\bfnm{P.}\binits{P.}} \AND
  \bauthor{\bsnm{Serfling},~\bfnm{R.}\binits{R.}}
(\byear{1988}).
\btitle{Strong uniform consistency rates for estimators of conditional
  functionals}.
\bjournal{Ann. Statist.}
\bvolume{16}
\bpages{1428--1449}.
\bid{doi={10.1214/aos/1176351047}, issn={0090-5364}, mr={0964932}}
\bptok{imsref}%
\end{barticle}
\endbibitem

\bibitem[\protect\citeauthoryear{Heckman, Ichimura and
  Todd}{1998}]{heckman1998matching}
\begin{barticle}[mr]
\bauthor{\bsnm{Heckman},~\bfnm{James~J.}\binits{J.~J.}},
  \bauthor{\bsnm{Ichimura},~\bfnm{Hidehiko}\binits{H.}} \AND
  \bauthor{\bsnm{Todd},~\bfnm{Petra}\binits{P.}}
(\byear{1998}).
\btitle{Matching as an econometric evaluation estimator}.
\bjournal{Rev. Econom. Stud.}
\bvolume{65}
\bpages{261--294}.
\bid{doi={10.1111/1467-937X.00044}, issn={0034-6527}, mr={1623713}}
\bptok{imsref}%
\end{barticle}
\endbibitem

\bibitem[\protect\citeauthoryear{Heckman and
  Vytlacil}{2005}]{heckman2005structural}
\begin{barticle}[mr]
\bauthor{\bsnm{Heckman},~\bfnm{James~J.}\binits{J.~J.}} \AND
  \bauthor{\bsnm{Vytlacil},~\bfnm{Edward}\binits{E.}}
(\byear{2005}).
\btitle{Structural equations, treatment effects, and econometric policy
  evaluation}.
\bjournal{Econometrica}
\bvolume{73}
\bpages{669--738}.
\bid{doi={10.1111/j.1468-0262.2005.00594.x}, issn={0012-9682}, mr={2135141}}
\bptok{imsref}%
\end{barticle}
\endbibitem

\bibitem[\protect\citeauthoryear{Imbens and Newey}{2009}]{imbens2009identification}
\begin{barticle}[mr]
\bauthor{\bsnm{Imbens},~\bfnm{Guido~W.}\binits{G.~W.}} \AND
  \bauthor{\bsnm{Newey},~\bfnm{Whitney~K.}\binits{W.~K.}}
(\byear{2009}).
\btitle{Identification and estimation of triangular simultaneous equations
  models without additivity}.
\bjournal{Econometrica}
\bvolume{77}
\bpages{1481--1512}.
\bid{doi={10.3982/ECTA7108}, issn={0012-9682}, mr={2561069}}
\bptok{imsref}%
\end{barticle}
\endbibitem

\bibitem[\protect\citeauthoryear{Kanaya and Kristensen}{2009}]{kristensen2009}
\begin{bmisc}[author]
\bauthor{\bsnm{Kanaya},~\bfnm{S.}\binits{S.}} \AND
  \bauthor{\bsnm{Kristensen},~\bfnm{D.}\binits{D.}}
(\byear{2009}).
\btitle{{{Estimation of stochastic volatility models by nonparametric
  filtering}}}.
\bhowpublished{Unpublished manuscript}.
\bptok{imsref}%
\end{bmisc}
\endbibitem

\bibitem[\protect\citeauthoryear{Lewbel and
  Linton}{2002}]{linton2002nonparametric}
\begin{barticle}[mr]
\bauthor{\bsnm{Lewbel},~\bfnm{Arthur}\binits{A.}} \AND
  \bauthor{\bsnm{Linton},~\bfnm{Oliver}\binits{O.}}
(\byear{2002}).
\btitle{Nonparametric censored and truncated regression}.
\bjournal{Econometrica}
\bvolume{70}
\bpages{765--779}.
\bid{doi={10.1111/1468-0262.00304}, issn={0012-9682}, mr={1913830}}
\bptok{imsref}%
\end{barticle}
\endbibitem

\bibitem[\protect\citeauthoryear{Li and
  Wooldridge}{2002}]{li2002semiparametric}
\begin{barticle}[mr]
\bauthor{\bsnm{Li},~\bfnm{Qi}\binits{Q.}} \AND
  \bauthor{\bsnm{Wooldridge},~\bfnm{Jeffrey~M.}\binits{J.~M.}}
(\byear{2002}).
\btitle{Semiparametric estimation of partially linear models for dependent data
  with generated regressors}.
\bjournal{Econometric Theory}
\bvolume{18}
\bpages{625--645}.
\bid{doi={10.1017/S0266466602183034}, issn={0266-4666}, mr={1906328}}
\bptok{imsref}%
\end{barticle}
\endbibitem

\bibitem[\protect\citeauthoryear{Linton and Nielsen}{1995}]{linton1995kernel}
\begin{barticle}[mr]
\bauthor{\bsnm{Linton},~\bfnm{Oliver}\binits{O.}} \AND
  \bauthor{\bsnm{Nielsen},~\bfnm{Jens~Perch}\binits{J.~P.}}
(\byear{1995}).
\btitle{A kernel method of estimating structured nonparametric regression based
  on marginal integration}.
\bjournal{Biometrika}
\bvolume{82}
\bpages{93--100}.
\bid{doi={10.1093/biomet/82.1.93}, issn={0006-3444}, mr={1332841}}
\bptok{imsref}%
\end{barticle}
\endbibitem

\bibitem[\protect\citeauthoryear{Mammen, Linton and
  Nielsen}{1999}]{mammen1999backfitting}
\begin{barticle}[mr]
\bauthor{\bsnm{Mammen},~\bfnm{E.}\binits{E.}},
  \bauthor{\bsnm{Linton},~\bfnm{O.}\binits{O.}} \AND
  \bauthor{\bsnm{Nielsen},~\bfnm{J.}\binits{J.}}
(\byear{1999}).
\btitle{The existence and asymptotic properties of a backfitting projection
  algorithm under weak conditions}.
\bjournal{Ann. Statist.}
\bvolume{27}
\bpages{1443--1490}.
\bid{doi={10.1214/aos/1017939137}, issn={0090-5364}, mr={1742496}}
\bptok{imsref}%
\end{barticle}
\endbibitem

\bibitem[\protect\citeauthoryear{Mammen, Rothe and
  Schienle}{2011}]{mammen2011semi}
\begin{bmisc}[author]
\bauthor{\bsnm{Mammen},~\bfnm{E.}\binits{E.}},
  \bauthor{\bsnm{Rothe},~\bfnm{C.}\binits{C.}} \AND
  \bauthor{\bsnm{Schienle},~\bfnm{M.}\binits{M.}}
(\byear{2011}).
\btitle{{Semiparametric estimation with generated covariates}}.
\bhowpublished{Unpublished manuscript}.
\bptok{imsref}%
\end{bmisc}
\endbibitem

\bibitem[\protect\citeauthoryear{Masry}{1996}]{masry1996multivariate}
\begin{barticle}[mr]
\bauthor{\bsnm{Masry},~\bfnm{Elias}\binits{E.}}
(\byear{1996}).
\btitle{Multivariate local polynomial regression for time series: Uniform
  strong consistency and rates}.
\bjournal{J. Time Ser. Anal.}
\bvolume{17}
\bpages{571--599}.
\bid{doi={10.1111/j.1467-9892.1996.tb00294.x}, issn={0143-9782}, mr={1424907}}
\bptok{imsref}%
\end{barticle}
\endbibitem

\bibitem[\protect\citeauthoryear{Newey}{1994a}]{newey1994kernel}
\begin{barticle}[mr]
\bauthor{\bsnm{Newey},~\bfnm{Whitney~K.}\binits{W.~K.}}
(\byear{1994}a).
\btitle{Kernel estimation of partial means and a general variance estimator}.
\bjournal{Econometric Theory}
\bvolume{10}
\bpages{233--253}.
\bid{doi={10.1017/S0266466600008409}, issn={0266-4666}, mr={1293201}}
\bptok{imsref}%
\end{barticle}
\endbibitem

\bibitem[\protect\citeauthoryear{Newey}{1994b}]{newey1994variance}
\begin{barticle}[mr]
\bauthor{\bsnm{Newey},~\bfnm{Whitney~K.}\binits{W.~K.}}
(\byear{1994}b).
\btitle{The asymptotic variance of semiparametric estimators}.
\bjournal{Econometrica}
\bvolume{62}
\bpages{1349--1382}.
\bid{doi={10.2307/2951752}, issn={0012-9682}, mr={1303237}}
\bptok{imsref}%
\end{barticle}
\endbibitem

\bibitem[\protect\citeauthoryear{Newey}{1997}]{newey1997convergence}
\begin{barticle}[mr]
\bauthor{\bsnm{Newey},~\bfnm{Whitney~K.}\binits{W.~K.}}
(\byear{1997}).
\btitle{Convergence rates and asymptotic normality for series estimators}.
\bjournal{J. Econometrics}
\bvolume{79}
\bpages{147--168}.
\bid{doi={10.1016/S0304-4076(97)00011-0}, issn={0304-4076}, mr={1457700}}
\bptok{imsref}%
\end{barticle}
\endbibitem

\bibitem[\protect\citeauthoryear{Newey, Powell and
  Vella}{1999}]{newey1999nonparametric}
\begin{barticle}[mr]
\bauthor{\bsnm{Newey},~\bfnm{Whitney~K.}\binits{W.~K.}},
  \bauthor{\bsnm{Powell},~\bfnm{James~L.}\binits{J.~L.}} \AND
  \bauthor{\bsnm{Vella},~\bfnm{Francis}\binits{F.}}
(\byear{1999}).
\btitle{Nonparametric estimation of triangular simultaneous equations models}.
\bjournal{Econometrica}
\bvolume{67}
\bpages{565--603}.
\bid{doi={10.1111/1468-0262.00037}, issn={0012-9682}, mr={1685723}}
\bptok{imsref}%
\end{barticle}
\endbibitem

\bibitem[\protect\citeauthoryear{Pagan}{1984}]{pagan1984econometric}
\begin{barticle}[mr]
\bauthor{\bsnm{Pagan},~\bfnm{Adrian}\binits{A.}}
(\byear{1984}).
\btitle{Econometric issues in the analysis of regressions with generated
  regressors}.
\bjournal{Internat. Econom. Rev.}
\bvolume{25}
\bpages{221--247}.
\bid{doi={10.2307/2648877}, issn={0020-6598}, mr={0741926}}
\bptok{imsref}%
\end{barticle}
\endbibitem

\bibitem[\protect\citeauthoryear{Song}{2008}]{song2008uniform}
\begin{barticle}[mr]
\bauthor{\bsnm{Song},~\bfnm{Kyungchul}\binits{K.}}
(\byear{2008}).
\btitle{Uniform convergence of series estimators over function spaces}.
\bjournal{Econometric Theory}
\bvolume{24}
\bpages{1463--1499}.
\bid{doi={10.1017/S0266466608080584}, issn={0266-4666}, mr={2456535}}
\bptok{imsref}%
\end{barticle}
\endbibitem

\bibitem[\protect\citeauthoryear{Sperlich}{2009}]{sperlich2009note}
\begin{barticle}[mr]
\bauthor{\bsnm{Sperlich},~\bfnm{Stefan}\binits{S.}}
(\byear{2009}).
\btitle{A note on non-parametric estimation with predicted variables}.
\bjournal{Econom. J.}
\bvolume{12}
\bpages{382--395}.
\bid{doi={10.1111/j.1368-423X.2009.00291.x}, issn={1368-4221}, mr={2562393}}
\bptok{imsref}%
\end{barticle}
\endbibitem

\bibitem[\protect\citeauthoryear{Stone}{1985}]{stone1985additive}
\begin{barticle}[mr]
\bauthor{\bsnm{Stone},~\bfnm{Charles~J.}\binits{C.~J.}}
(\byear{1985}).
\btitle{Additive regression and other nonparametric models}.
\bjournal{Ann. Statist.}
\bvolume{13}
\bpages{689--705}.
\bid{doi={10.1214/aos/1176349548}, issn={0090-5364}, mr={0790566}}
\bptok{imsref}%
\end{barticle}
\endbibitem

\bibitem[\protect\citeauthoryear{van~de Geer}{2000}]{vandegeer2009book}
\begin{bbook}[author]
\bauthor{\bparticle{van~de} \bsnm{Geer},~\bfnm{S.}\binits{S.}}
(\byear{2000}).
\btitle{Empirical Processes in M-Estimation}.
\bpublisher{Cambridge Univ. Press}, \baddress{Cambridge}.
\bptok{imsref}%
\end{bbook}
\endbibitem

\bibitem[\protect\citeauthoryear{van~der Vaart and Wellner}{1996}]{van1996weak}
\begin{bbook}[mr]
\bauthor{\bparticle{van~der} \bsnm{Vaart},~\bfnm{Aad~W.}\binits{A.~W.}} \AND
  \bauthor{\bsnm{Wellner},~\bfnm{Jon~A.}\binits{J.~A.}}
(\byear{1996}).
\btitle{Weak Convergence and Empirical Processes: With Applications to Statistics}.
\bpublisher{Springer}, \baddress{New York}.
\bid{mr={1385671}}
\bptok{imsref}%
\end{bbook}
\endbibitem

\end{thebibliography}
\end{document}